\documentclass[12pt,reqno]{amsart}

\usepackage[arrow,matrix,curve]{xy}

\usepackage[dvips]{graphicx} 

\usepackage{amssymb, latexsym, amsmath, amscd, array, 
}

\newtheorem{theorem}{Theorem}[section]

\theoremstyle{definition}

\newtheorem{question}[theorem]{Question}

\numberwithin{equation}{section}

\newcommand\N {{\mathbb N}} 

\newcommand\R {{\mathbb R}}
\newcommand\Q {{\mathbb Q}}
\newcommand\Z {{\mathbb Z}}

\newcommand\RRR{{\mbox{I\!I\!R}}}

\begin{document}

\thispagestyle{empty}

\title[A Burgessian critique of nominalistic tendencies]{A Burgessian
critique of nominalistic tendencies in contemporary mathematics and
its historiography}

\author[K. Katz and M. Katz]{Karin Usadi Katz and Mikhail G. Katz$^{0}$}

\address{Department of Mathematics, Bar Ilan University, Ramat Gan
52900 Israel} \email{katzmik@macs.biu.ac.il}

\footnotetext{Supported by the Israel Science Foundation grant
1294/06}



\bigskip\bigskip\bigskip\noindent

\begin{abstract}
We analyze the developments in mathematical rigor from the viewpoint
of a Burgessian critique of nominalistic reconstructions.  We apply
such a critique to the reconstruction of infinitesimal analysis
accomplished through the efforts of Cantor, Dedekind, and Weierstrass;
to the reconstruction of Cauchy's foundational work associated with
the work of Boyer and Grabiner; and to Bishop's constructivist
reconstruction of classical analysis.  We examine the effects of an
ontologically limitative disposition on historiography, teaching, and
research.
\end{abstract}

\maketitle

\tableofcontents

\section{Introduction}

Over the course of the past 140 years, the field of professional pure
mathematics (analysis in particular), and to a large extent also its
professional historiography, have become increasingly dominated by a
particular philosophical disposition.  We argue that such a
disposition is akin to nominalism, and  examine its ramifications.

In 1983, J.~Burgess proposed a useful dichotomy for analyzing
nominalistic narratives.  The starting point of his critique is his
perception that a philosopher's job is not to rule on the ontological
merits of this or that scientific entity, but rather to try to
understand those entities that are employed in our best scientific
theories.  From this viewpoint, the problem of nominalism is the
awkwardness of the contortions a nominalist goes through in developing
an alternative to his target scientific practice, an alternative
deemed ontologically ``better" from his reductive perspective, but in
reality amounting to the imposition of artificial strictures on the
scientific practice.

Burgess introduces a dichotomy of {\em hermeneutic\/} versus {\em
revolutionary\/} nominalism.  Thus, {\em hermeneutic nominalism\/} is
the hypothesis that science, properly interpreted, already dispenses
with mathematical objects (entities) such as numbers and sets.
Meanwhile, {\em revolutionary nominalism\/} is the project of
replacing current scientific theories by alternatives dispensing with
mathematical objects, see Burgess \cite[p.~96]{Bu83} and Burgess and
Rosen~\cite{BR}.

Nominalism in the philosophy of mathematics is often understood
narrowly, as exemplified by the ideas of J. S. Mill and P. Kitcher,
going back to Aristotle.%
\footnote{See for example S.~Shapiro \cite{Sha}.}
However, the Burgessian distinction between hermeneutic and
revolutionary reconstructions can be applied more broadly, so as to
include nominalistic-type reconstructions that vary widely in their
ontological target, namely the variety of abstract objects (entities)
they seek to challenge (and, if possible, eliminate) as being merely
conventional, see \cite[p.~98-99]{Bu83}.

Burgess quotes at length Yu.~Manin's critique in the 1970s of
mathematical nominalism of constructivist inspiration, whose
ontological target is the classical infinity, namely, 
\begin{quote}
abstractions which are infinite and do not lend themselves to a
constructivist interpretation \cite[p.~172-173]{Man}.
\end{quote}
This suggests that Burgess would countenance an application of his
dichotomy to nominalism of a constructivist inspiration.

The ontological target of the constructivists is the concept of
Cantorian infinities, or more fundamentally, the logical principle of
the Law of Excluded Middle (LEM).  Coupled with a classical
interpretation of the existence quantifier, LEM is responsible for
propelling the said infinities into a dubious existence.  LEM is the
abstract object targeted by Bishop's constructivist nominalism, which
can therefore be called an anti-LEM nominalism.%
\footnote{A more detailed discussion of LEM in the context of the
proof of the irrationality of~$\sqrt{2}$ may be found in
Section~\ref{two}, see footnote~\ref{root}.}
Thus, anti-LEM nominalism falls within the scope of the Burgessian
critique, and is the first of the nominalistic reconstructions we wish
to analyze.

The anti-LEM nominalistic reconstruction was in fact a {\em
re\/}-recon\-struction of an earlier nominalistic reconstruction of
analysis, dating from the 1870s.  The earlier reconstruction was
implemented by the great triumvirate%
\footnote{\label{triumvirate}C.~Boyer refers to Cantor, Dedekind, and
Weierstrass as ``the great triumvirate'', see \cite[p.~298]{Boy}.}
of Cantor, Dedekind, and Weierstrass.  The ontological target of the
triumvirate reconstruction was the abstract entity called the {\em
infinitesimal\/}, a basic building block of a continuum,%
\footnote{For an entertaining history of infinitesimals see P.~Davis
and R.~Hersh \cite[p.~237-254]{DH}.  For an analysis of a variety of
competing theories of the continuum, see P.~Ehrlich~\cite{Eh94} as
well as R.~Taylor~\cite{Ta}.}
according to a line of investigators harking back to the Greek
antiquity.%
\footnote{\label{klein}See also footnote~\ref{cauchy} on Cauchy.}

To place these historical developments in context, it is instructive
to examine Felix Klein's remarks dating from 1908.  Having outlined
the developments in real analysis associated with Weierstrass and his
followers, Klein pointed out that 
\begin{quote}
The scientific mathematics of today is built upon the series of
developments which we have been outlining.  But an essentially
different conception of infinitesimal calculus has been running
parallel with this [conception] through the centuries
\cite[p.~214]{Kl}.
\end{quote}
Such a different conception, according to Klein,
\begin{quote}
harks back to old metaphysical speculations concerning the {\em
structure of the continuum\/} according to which this was made up of
[...] infinitely small parts \cite[p.~214]{Kl} [emphasis
added---authors].
\end{quote}
The significance of the triumvirate reconstruction has often been
measured by the yardstick of the extirpation of the infinitesimal.%
\footnote{\label{hob}Thus, after describing the formalisation of the
real continuum in the 1870s, on pages 127-128 of his retiring
presidential address in 1902, E.~Hobson remarks triumphantly as
follows: ``It should be observed that the criterion for the
convergence of an aggregate [i.e. an equivalence class defining a real
number] is of such a character that no use is made in it of {\em
infinitesimals\/}'' \cite[p.~128]{Ho} [emphasis added--authors].
Hobson reiterates: ``In all such proofs [of convergence] the only
statements made are as to relations of finite numbers, no such
entities as {\em infinitesimals\/} being recognized or employed.  Such
is the essence of the~$\epsilon$[$, \delta$] proofs with which we are
familiar'' \cite[p.~128]{Ho} [emphasis added--authors].  The tenor of
Hobson's remarks is that Weierstrass's fundamental accomplishment was
the elimination of infinitesimals from foundational discourse in
analysis.}

The infinitesimal ontological target has similarly been the motivating
force behind a more recent nominalistic reconstruction, namely a
nominalistic re-appraisal of the meaning of Cauchy's foundational work
in analysis.

We will analyze these three nominalistic projects through the lens of
the dichotomy introduced by Burgess.  Our preliminary conclusion is
that, while the triumvirate reconstruction was primarily revolutionary
in the sense of Burgess, and the (currently prevailing) Cauchy
reconstruction is mainly hermeneutic, the anti-LEM reconstruction has
combined elements of both types of nominalism.  We will examine the
effects of a nominalist disposition on historiography, teaching, and
research.

A traditional view of 19th century analysis holds that a search for
rigor inevitably leads to epsilontics, as developed by Weierstrass in
the 1870s; that such inevitable developments culminated in the
establishment of ultimate set-theoretic foundations for mathematics by
Cantor; and that eventually, once the antinomies sorted out, such
foundations were explicitly expressed in axiomatic form by Zermelo and
Fraenkel.  Such a view entails a commitment to a specific destination
or ultimate goal of scientific devepment as being pre-determined and
intrinsically inevitable.  The postulation of a specific outcome,
believed to be the inevitable result of the development of the
discipline, is an outlook specific to the mathematical community.
Challenging such a {\em belief\/} appears to be a radical proposition
in the eyes of a typical professional mathematician, but not in the
eyes of scientists in related fields of the exact sciences.  It is
therefore puzzling that such a view should be accepted without
challenge by a majority of historians of mathematics, who tend to toe
the line on the mathematicians' belief.  Could mathematical analysis
have followed a different path of development?  Related material
appears in Alexander \cite{Al}, Giordano \cite{Gio}, Katz and Tall
\cite{KT}, Kutateladze \cite[chapter~63]{Ku}, Mormann \cite{Mo},
Sepkoski \cite{Se}, and Wilson \cite{Wi}.

\section{An anti-LEM nominalistic reconstruction}
\label{glossary}
\label{two}

This section is concerned with E.~Bishop's approach to reconstructing
analysis.  Bishop's approach is rooted in Brouwer's revolt against the
non-constructive nature of mathematics as practiced by his
contemporaries.%
\footnote{Similar tendencies on the part of Wittgenstein were analyzed
by H.~Putnam, who describes them as ``minimalist''
\cite[p.~242]{Pu07}.  See also G.~Kreisel \cite{Kr}.}

Is there meaning after LEM?  The Brouwer--Hilbert debate captured the
popular mathematical imagination in the 1920s.  Brouwer's crying call
was for the elimination of most of the applications of LEM from
meaningful mathematical discourse.  Burgess discusses the debate
briefly in his treatment of nominalism in \cite[p.~27]{Bu04}.  We will
analyze E.~Bishop's implementation of Brouwer's nominalistic project.%
\footnote{It has been claimed that Bishopian constructivism, unlike
Brouwer's intuitionism, is compatible with classical mathematics, see
e.g.~Davies \cite{Da05}.  However, Brouwerian counterexamples do
appear in Bishop's work; see footnote~\ref{lpo} for more details.
This could not be otherwise, since a verificational interpretation of
the quantifiers necessarily results in a clash with classical
mathematics.  As a matter of presentation, the conflict with classical
mathematics had been de-emphasized by Bishop.  Bishop finesses the
issue of Brouwer's theorems (e.g., that every function is continuous)
by declaring that he will only deal with uniformly continuous
functions to begin with.  In Bishopian mathematics, a circle cannot be
decomposed into a pair of antipodal sets.  A counterexample to the
classical extreme value theorem is discussed in \cite[p.~295]{TV}, see
footnote~\ref{lpo} for details.}

It is an open secret that the much-touted success of Bishop's
implementation of the intuitionistic project in his 1967 book
\cite{Bi67} is due to philosophical compromises with a Platonist
viewpoint that are resolutely rejected by the intuitionistic
philosopher M.~Dummett \cite{Du77}.  Thus, in a dramatic departure
from both Kronecker%
\footnote{Kronecker's position is discussed in more detail in
Section~\ref{hersh}, in the main text around footnote~\ref{boniface}.}
and Brouwer, Bishopian constructivism accepts the completed (actual)
infinity of the integers~$\Z$.%
\footnote{\label{potential}Intuitionists view~$\N$ as a potential
totality; for a more detailed discussion see, e.g.,
\cite[section~4.3]{Fef}.}

Bishop expressed himself as follows on the first page of his book:
\begin{quote}
in another universe, with another biology and another physics, [there]
will develop mathematics which in essence is the same as ours
\cite[p.~1]{Bi67}.  
\end{quote}
Since the sensory perceptions of the human {\em body\/} are physics-
and chemistry-bound, a claim of such trans-universe invariance%
\footnote{a claim that we attribute to a post-{\em Sputnik\/} fever}
amounts to the positing of a {\em disembodied\/} nature of the
infinite natural number system, transcending physics and chemistry.%
\footnote{The muddle of realism and anti-realism in the anti-LEM
sector will be discussed in more detail in Section~\ref{insurrection}.
Bishop's disembodied integers illustrate the awkward philosophical
contorsions which are a tell-tale sign of nominalism.  An alternative
approach to the problem is pursued in modern cognitive science.
Bishop's disembodied integers, the cornerstone of his approach, appear
to be at odds with modern cognitive theory of {\em embodied\/}
knowledge, see Tall~\cite{Ta09a}, Lakoff and N\'u\~nez~\cite{LN},
Sfard~\cite{Sf}, Yablo~\cite{Ya}.  Reyes \cite{Re} presents an
intriguing thesis concerning an allegedly rhetorical nature of
Newton's attempts at grounding infinitesimals in terms of {\em
moments\/} or {\em nascent and evanescent quantities\/}, and Leibniz's
similar attempts in terms of ``a general heuristic under [the name of]
the principle of continuity'' \cite[p.~172]{Re}.  He argues that what
made these theories vulnerable to criticism is the reigning principle
in 17th century methodology according to which abstract objects must
necessarily have empirical counterparts/referents.  D.~Sherry points
out that ``Formal axiomatics emerged only in the 19th century, after
geometry embraced objects with no empirical counterparts (e.g.,
Poncelet's points at infinity...)''  \cite[p.~67]{She09}.  See also
S.~Feferman's approach of conceptual structuralism \cite{Fef}, for a
view of mathematical objects as mental conceptions.}

What type of nominalistic reconstruction best fits the bill of
Bishop's constructivism?  Bishop's rejection of classical mathematics
as a ``debasement of meaning''%
\footnote{Bishop diagnosed classical mathematics with a case of a
``debasement of meaning'' in his {\em Schizophrenia in contemporary
mathematics\/} (1973).  Hot on the heels of {\em Schizophrenia\/} came
the 1975 {\em Crisis in contemporary mathematics\/} \cite{Bi75}, where
the same diagnosis was slapped upon infinitesimal calculus \`a la
Robinson \cite{Ro66}.}
\cite[p.~1]{Bi85} would place him squarely in the camp of
revolutionary nominalisms in the sense of Burgess; yet some elements
of Bishop's program tend to be of the hermeneutic variety, as well.

As an elementary example, consider Bishop's discussion of the
irrationality of the square root of~$2$ in \cite[p.~18]{Bi85}.
Irrationality is defined constructively in terms of being quantifiably
apart from each rational number.  The classical proof of the
irrationality of~$\sqrt{2}$ is a proof by contradiction.  Namely, we
{\em assume\/} a hypothesized equality~$\sqrt{2}=\frac{m}{n}$, examine
the parity of the powers of~$2$, and arrive at a contradiction.  At
this stage, irrationality is considered to have been proved, in
classical logic.%
\footnote{\label{f1}The classical proof showing that~$\sqrt{2}$ is
{\em not rational\/} is, of course, acceptable in intuitionistic
logic.  To pass from this to the claim of its {\em irrationality\/} as
defined above, requires LEM (see footnote~\ref{root} for details).}

However, as Bishop points out, the proof can be modified slightly so
as to avoid LEM, and acquire an enhanced {\em numerical meaning\/}.
Thus, {\em without\/} exploiting the equality~$\sqrt{2}=\frac{m}{n}$,
one can exhibit effective positive lower bounds for the
difference~$|\sqrt{2}-\frac{m}{n}|$ in terms of the denominator~$n$,
resulting in a constructively adequate proof of irrationality.%
\footnote{\label{root}Such a proof may be given as follows.  For each
rational~$m/n$, the integer~$2n^2$ is divisible by an odd power
of~$2$, while~$m^2$ is divisible by an even power of~$2$.
Hence~$|2n^2-m^2|\geq 1$ (here we have applied LEM to an effectively
decidable predicate over~$\Z$, or more precisely the law of
trichotomy).  Since the decimal expansion of~$\sqrt{2}$ starts
with~$1.41\ldots$, we may assume~$\frac{m}{n} \leq 1.5$.  It follows
that
\begin{equation*}
|\sqrt{2} - \tfrac{m}{n}| = \frac{|2n^2 - m^2|}{n^2 \left(
  \sqrt{2}+\tfrac{m}{n} \right)} \geq \frac{1}{n^2
  \left(\sqrt{2}+\tfrac{m}{n}\right)} \geq \frac{1}{3n^2},
\end{equation*}
yielding a numerically meaningful proof of irrationality, which is a
special case of Liouville's theorem on diophantine approximation of
algebraic numbers, see \cite{HW}.}
Such a proof is merely a modification of a classical proof, and can
thus be considered a hermeneutic reconstruction thereof.  A number of
classical results (though by no means all) can be reinterpreted
constructively, resulting in an enhancement of their numerical
meaning, in some cases at little additional cost.  This type of
project is consistent with the idea of a hermeneutic nominalism in the
sense of Burgess, and related to the notion of {\em liberal
constructivism\/} in the sense of G.~Hellman (see below).

The intuitionist/constructivist opposition to classical mathematics is
predicated on the the philosophical assumption that ``meaningful"
mathematics is mathematics done without the law of excluded middle.
E.~Bishop (following Brouwer but surpassing him in rhetoric) is on
record making statements of the sort 
\begin{quote}
``Very possibly classical mathematics will cease to exist as an
independent discipline'' \cite[p.~54]{Bi68}
\end{quote}
(to be replaced, naturally, by constructive mathematics); and
\begin{quote}
``Brouwer's criticisms of classical mathematics were concerned with
what I shall refer to as `the debasement of meaning' ''
\cite[p.~1]{Bi85}.
\end{quote}
Such a stance posits intuitionism/constructivism as an {\em
alternative\/} to classical mathematics, and is described as {\em
radical constructivism\/} by G.~Hellman~\cite[p.~222]{Hel93a}.
Radicalism is contrasted by Hellman with a {\em liberal\/} brand of
intuitionism (a {\em companion\/} to classical mathematics).  Liberal
constructivism may be exemplified by A.~Heyting \cite{He, He73}, who
was Brouwer's student, and formalized intuitionistic logic.

To motivate the long march through the foundations occasioned by a
LEM-eliminative agenda, Bishop \cite{Bi75} goes to great lengths to
dress it up in an appealing package of a {\em theory of meaning\/}
that first conflates meaning with numerical meaning (a goal many
mathematicians can relate to), and then numerical meaning with LEM
extirpation.%
\footnote{Such a reduction is discussed in more detail in
Section~\ref{hersh}.}
Rather than merely rejecting LEM or related logical principles such as
trichotomy which sound perfectly unexceptionable to a typical
mathematician, Bishop presents these principles in quasi metaphysical
garb of ``principles of omniscience".%
\footnote{\label{lpo}Thus, the main target of his criticism in
\cite{Bi75} is the ``limited principle of omniscience'' (LPO).  The
LPO is formulated in terms of sequences, as the principle that it is
possible to search ``a sequence of integers to see whether they all
vanish'' \cite[p.~511]{Bi75}.  The LPO is equivalent to the law of
trichotomy:~$(a< 0) \vee (a=0) \vee (a> 0)$.  An even weaker principle
is~$(a\leq 0) \vee (a \geq 0)$, whose failure is exploited in the
construction of a counterexample to the extreme value theorem by
Troelstra and van Dalen \cite[p.~295]{TV}, see also our
Section~\ref{hersh}.  This property is false intuitionistically.
After discussing real numbers~$x\geq 0$ such that it is ``{\em
not\/}'' true that~$x>0$ or~$x=0$, Bishop writes:
\begin{quote}
In much the same way we can construct a real number~$x$ such that it
is {\em not\/} true that~$x\geq 0$ or~$x\leq 0$ \cite[p.~26]{Bi67},
\cite[p.~28]{BB}.
\end{quote}
An~$a$ satisfying~$\neg ((a\leq 0) \vee (a \geq 0))$ immediately
yields a counterexample~$f(x)=ax$ to the extreme value theorem (EVT)
on~$[0,1]$ (see \cite[p.~59, exercise~9]{Bi67}; \cite[p.~62,
exercise~11]{BB}).  Bridges interprets Bishop's italicized ``{\em
not\/}'' as referring to a Brouwerian counterexample, and asserts that
trichotomy as well as the principle~$(a\leq 0) \vee (a \geq 0)$ are
independent of Bishopian constructivism.  See D.~Bridges~\cite{Br94}
for details; a useful summary may be found in Taylor \cite{Ta}.}
Bishop retells a creation story of intuitionism in the form of an
imaginary dialog between Brouwer and Hilbert where the former
completely dominates the exchange.  Indeed, Bishop's imaginary
Brouwer-Hilbert exchange is dominated by an unspoken assumption that
Brouwer is the only one who seeks ``meaning", an assumption that his
illustrious opponent is never given a chance to challenge.  Meanwhile,
Hilbert's comments in 1919 reveal clearly his attachment to meaning
which he refers to as {\em internal necessity\/}:
\begin{quote}
We are not speaking here of arbitrariness in any sense.  Mathematics
is not like a game whose tasks are determined by arbitrarily
stipulated rules. Rather, it is a conceptual system possessing
internal necessity that can only be so and by no means otherwise
\cite[p.~14]{Hi} (cited in Corry \cite{Cor}).
\end{quote}

A majority of mathematicians (including those favorable to
constructivism) feel that an implementation of Bishop's program does
involve a significant complication of the technical development of
analysis, as a result of the nominalist work of LEM-elimination.
Bishop's program has met with a certain amount of success, and
attracted a number of followers.  Part of the attraction stems from a
detailed lexicon developed by Bishop so as to challenge received
(classical) views on the nature of mathematics.  A constructive
lexicon was a {\em sine qua non\/} of his success.  A number of terms
from Bishop's constructivist lexicon constitute a novelty as far as
intuitionism is concerned, and are not necessarily familiar even to
someone knowledgeable about intuitionism {\em per se\/}.  It may be
helpful to provide a summary of such terms for easy reference,
arranged alphabetically, as follows.

\medskip
$\bullet$ {\em Debasement of meaning\/} is the cardinal sin of the
classical opposition, from Cantor to Keisler,%
\footnote{\label{keisler1}But see footnote~\ref{keisler2}.}
committed with {\em LEM\/} (see below).  The term occurs in Bishop's
{\em Schizophrenia\/} \cite{Bi85} and {\em Crisis\/}~\cite{Bi75}
texts.

\medskip
$\bullet$ {\em Fundamentalist excluded thirdist\/} is a term that
refers to a classically-trained mathematician who has not yet become
sensitized to implicit use of the law of excluded middle (i.e.,~{\em
excluded third\/}) in his arguments, see \cite[p.~249]{Ri}.%
\footnote{\label{rich1}This use of the term ``fundamentalist excluded
thirdist'' is in a text by Richman, not Bishop.  I have not been able
to source its occurrence in Bishop's writing.  In a similar vein, an
ultrafinitist recently described this writer as a ``choirboy of
infinitesimology''; however, this term does not seem to be in general
use.  See also footnote~\ref{rich2}.}

\medskip
$\bullet$ {\em Idealistic mathematics\/} is the output of Platonist
mathematical sensibilities, abetted by a metaphysical faith in {\em
LEM\/} (see below), and characterized by the presence of merely a {\em
peculiar pragmatic content\/} (see below).

\medskip
$\bullet$ {\em Integer\/} is the revealed source of all {\em
meaning\/} (see below), posited as an alternative foundation
displacing both formal logic, axiomatic set theory, and recursive
function theory.  The integers wondrously escape%
\footnote{By dint of a familiar oracular quotation from Kronecker; see
also main text around footnote~\ref{potential}.}
the vigilant scrutiny of a constructivist intelligence determined to
uproot and nip in the bud each and every Platonist fancy of a concept
{\em external\/} to the mathematical mind.

\medskip
$\bullet$ {\em Integrity\/} is perhaps one of the most misunderstood
terms in Errett Bishop's lexicon.  Pourciau in his {\em Education\/}
\cite{Po99} appears to interpret it as an indictment of the ethics of
the classical opposition.  Yet in his {\em Schizophrenia\/} text,
Bishop merely muses:
\begin{quote}
[...] I keep coming back to the term ``integrity''.  \cite[p.~4]{Bi85}
\end{quote}
Note that the period is in the original.  Bishop describes {\em
integrity\/} as the opposite of a syndrome he colorfully refers to as
{\em schizophrenia\/}, characterized 
\begin{enumerate}
\item[(a)] by a rejection of common sense in favor of formalism, 
\item[(b)]
by {\em debasement of meaning\/} (see above), 
\item[(c)] 
as well as by a list of other ills---
\end{enumerate}
but {\em excluding\/} dishonesty.  Now the root of
\begin{equation*}
\mbox{integr-ity}
\end{equation*}
is identical with that of {\em integer\/} (see above), the Bishopian
ultimate foundation of analysis.  Bishop's evocation of {\em
integrity\/} may have been an innocent pun intended to allude to a
healthy constructivist mindset, where the {\em integers\/} are
uppermost.%
\footnote{In Bishop's system, the integers are uppermost to the
exclusion of the continuum.  Bishop rejected Brouwer's work on an
intuitionstic continuum in the following terms: 
\begin{quote}
Brouwer's bugaboo has been compulsive speculation about the nature of
the continuum.  His fear seems to have been that, unless he personally
intervened to prevent it, the continuum would turn out to be discrete.
[The result was Brouwer's] semimystical theory of the continuum
\cite[p.~6 and 10]{Bi67}.
\end{quote} 
Brouwer sought to incorporate a theory of the continuum as part of
intuitionistic mathematics, by means of his free choice sequences.
Bishop's commitment to integr-ity is thus a departure from Brouwerian
intuitionism.}

\medskip
$\bullet$ {\em Law of excluded middle (LEM)\/} is the main source of
the non-constructivities of classical mathematics.%
\footnote{See footnote~\ref{root} and footnote~\ref{lpo} for some
examples.}
Every formalisation of {\em intuitionistic logic\/} excludes {\em
LEM\/}; adding {\em LEM\/} back again returns us to {\em classical
logic\/}.

\medskip
$\bullet$ {\em Limited principle of omniscience (LPO)\/} is a weak
form of {\em LEM\/} (see above), involving {\em LEM\/}-like oracular
abilities limited to the context of integer sequences.%
\footnote{See footnote~\ref{lpo} for a discussion of LPO.}
The {\em LPO\/} is still unacceptable to a constructivist, but could
have served as a basis for a {\em meaningful\/} dialog between Brouwer
and Hilbert (see \cite{Bi75}), that could allegedly have changed the
course of 20th century mathematics.

\medskip
$\bullet$ {\em Meaning\/} is a favorite philosophic term in Bishop's
lexicon, necessarily preceding an investigation of {\em truth\/} in
any coherent discussion.  In Bishop's writing, the term {\em
meaning\/} is routinely conflated with {\em numerical meaning\/} (see
below).

\medskip
$\bullet$ {\em Numerical meaning\/} is the content of a theorem
admitting a proof based on intuitionistic logic, and expressing
computationally meaningful facts about the integers.%
\footnote{As an illustration, a numerically meaningful proof of the
irrationality of~$\sqrt{2}$ appears in footnote~\ref{root}.}
The conflation of {\em numerical meaning\/} with {\em meaning\/} par
excellence in Bishop's writing, has the following two consequences:
\begin{enumerate}
\item
it empowers the constructivist to sweep under the rug the distinction
between pre-LEM and post-LEM numerical meaning, lending a marginal
degree of plausibility to a dismissal of classical theorems which
otherwise appear eminently coherent and meaningful;%
\footnote{See the main text around footnote~\ref{boniface} for a
discussion of the classical extreme value theorem and its LEMless
remains.}
and
\item
it allows the constructivist to enlist the support of {\em
anti-realist\/} philosophical schools of thought (e.g.~Michael
Dummett) in the theory of meaning, inspite of the apparent tension
with Bishop's otherwise {\em realist\/} declarations (see entry {\em
realistic mathematics\/} below).
\end{enumerate}

\medskip
$\bullet$ {\em Peculiar pragmatic content\/} is an expression of
Bishop's \cite[p.~viii]{Bi67} that was analyzed by Billinge
\cite[p.~179]{Bil03}.  It connotes an alleged lack of empirical
validity of classical mathematics, when classical results are merely
{\em inference tickets\/} \cite[p.~180]{Bil03} used in the deduction
of other mathematical results.

\medskip
$\bullet$ {\em Realistic mathematics\/}.  The dichotomy of ``realist''
{\em versus\/} ``idealist'' (see above) is the dichotomy of
``constructive'' {\em versus\/} ``classical'' mathematics, in Bishop's
lexicon.  There are two main narratives of the Intuitionist
insurrection, one {\em anti-realist\/} and one {\em realist\/}.  The
issue is discussed in the next section.

\section{Insurrection according to Errett and according to Michael}
\label{insurrection}

The {\em anti-realist\/} narrative, mainly following Michael Dummett
\cite{Du}, traces the original sin of classical mathematics with {\em
LEM\/}, all the way back to Aristotle.%
\footnote{\label{keisler2}the entry under {\em debasement of
meaning\/} in Section~\ref{glossary} would read, accordingly, ``the
classical opposition from Aristotle to Keisler''; see main text at
footnote~\ref{keisler1}.}
The law of excluded middle (see Section~\ref{glossary}) is the
mathematical counterpart of geocentric cosmology (alternatively, of
phlogiston, see \cite[p.~299]{Po}), slated for the dustbin of
history.%
\footnote{\label{dummett}Following Kronecker and Brouwer, Dummett
rejects actual infinity, at variance with Bishop.}
The anti-realist narrative dismisses the Quine-Putnam indispensability
thesis (see Feferman \cite[Section IIB]{Fe}) on the grounds that a
{\em philosophy-first\/} examination of first principles is the unique
authority empowered to determine the correct way of doing
mathematics.%
\footnote{\label{doom}In Hellman's view, ``any [...]  attempt to
reinstate a `first philosophical' theory of meaning prior to all
science is doomed'' \cite[p.~439]{He98}.  What this appears to mean is
that, while there can certainly be a philosophical notion of meaning
before science, any attempt to {\em prescribe\/} standards of meaning
{\em prior\/} to the actual practice of science, is {\em doomed\/}.}
Generally speaking, it is this narrative that seems to be favored by a
number of philosophers of mathematics.

Dummett opposes a truth-valued, bivalent semantics, namely the notion
that truth is one thing and knowability another, on the grounds that
it violates Dummett's {\em manifestation requirement\/}, see
Shapiro~\cite[p.~54]{Sha01}.  The latter requirement, in the context
of mathematics, is merely a {\em restatement\/} of the intuitionistic
principle that truth is tantamount to verifiability (necessitating a
constructive interpretation of the quantifiers).  Thus, an acceptance
of Dummett's manifestation requirement, leads to intuitionistic
semantics and a rejection of LEM.

In his 1977 foundational text \cite{Du77} originating from 1973
lecture notes, Dummett is frank about the source of his interest in
the intuitionist/classical dispute in mathematics \cite[p.~ix]{Du77}:
\begin{quote} 
This dispute bears a {\bf strong resemblance} to other disputes over
realism of one kind or another, that is, concerning various kinds of
subject-matter (or types of statement), including that over realism
about the physical universe [emphasis added--authors]
\end{quote}
What Dummett proceeds to say at this point, reveals the nature of his
interest:
\begin{quote} 
but intuitionism represents the only sustained attempt by the
opponents of a {\bf realist view} to work out a coherent embodiment of
their philosophical beliefs [emphasis added--authors]
\end{quote}
What interests Dummett here is his fight against the {\em realist
view\/}.  What endears intuitionists to him, is the fact that they
have succeeded where the phenomenalists have not \cite[p.~ix]{Du77}:
\begin{quote}
Phenomenalists might have attained a greater success if they had made
a remotely comparable effort to show in detail what consequences their
interpretation of material-object statements would have for our
employment of our language.
\end{quote}
However, Dummett's conflation of the mathematical debate and the
philosophical debate, could be challenged.

We hereby explicitly sidestep the debate opposing the realist (as
opposed to the super-realist, see W.~Tait \cite{Ta01}) position and
the anti-realist position.  On the other hand, we observe that a
defense of indispensability of mathematics would necessarily start by
challenging Dummett's ``manifestation''.  More precisely, such a
defense would have to start by challenging the extension of Dummett's
manifestation requirement, from the realm of philosophy to the realm
of mathematics.  While Dummett chooses to pin the opposition to
intuitionism, to a belief \cite[p.~ix]{Du77} in an
\begin{quote}
interpretation of mathematical statements as referring to an
independently existing and objective reality[,]
\end{quote}
(i.e. a Platonic world of mathematical entities),
J.~Avigad~\cite{Av07} memorably retorts as follows:
\begin{quote}
We do not need fairy tales about numbers and triangles prancing about
in the realm of the abstracta.%
\footnote{\label{rich2}Constructivist Richman takes a dimmer view of
prancing numbers and triangles.  In addition, he presents a proposal
to eliminate the axiom of choice altogether from constructive
mathematics (including countable choice).  Since the ultrafilter axiom
is weaker than the axiom of choice, one might have hoped it would be
salvaged; not so:
\begin{quote}
We are all Platonists, aren't we? In the trenches, I mean---when the
chips are down. Yes, Virginia, there really are circles, triangles,
numbers, continuous functions, and all the rest.  Well, maybe not free
ultrafilters. Is it important to believe in the existence of free
ultrafilters? Surely that's not required of a Platonist. I can more
easily imagine it as a test of sanity: `He believes in free
ultrafilters, but he seems harmless' (Richman \cite{Ri94}).
\end{quote}
For Richman's contribution to constructivist lexicon see
footnote~\ref{rich1}.}
\end{quote}

Meanwhile, the {\em realist\/} narrative of the intuitionist
insurrection appears to be more consistent with what Bishop himself
actually wrote.  In his foundational essay, Bishop expresses his
position as follows:
\begin{quote}
As pure mathematicians, we must decide whether we are playing a game,
or whether our theorems describe an external reality
\cite[p.~507]{Bi75}.
\end{quote}
The right answer, to Bishop, is that they do describe an external
reality.%
\footnote{Our purpose here is not to endorse or refute Bishop's views
on this point, but rather to document his actual position, which
appears to diverge from Dummett's.}
The dichotomy of ``realist'' {\em versus\/} ``idealist'' is the
dichotomy of ``constructive'' {\em versus\/} ``classical''
mathematics, in Bishop's lexicon (see entry under {\em idealistic
mathematics\/} in Section~\ref{glossary}).  Bishop's ambition is to
incorporate ``such mathematically oriented disciplines as physics''
\cite[p.~4]{Bi85} as part of his constructive revolution, revealing a
recognition, on his part, of the potency of the Quine-Putnam
indispensability challenge.%
\footnote{Billinge \cite[p.~314]{Bil} purports to detect ``inchoate''
anti-realist views in Bishop's writings, but provides no constructive
proof of their existence, other than a pair quotes on {\em numerical
meaning\/}.  Meanwhile, Hellman \cite[p.~222]{Hel93a} writes: ``Some
of Bishop's remarks (1967) suggest that his position belongs in [the
radical] category''.}

N.~Kopell and G.~Stolzenberg, close associates of Bishop, published a
three-page {\em Commentary\/} \cite{KS} following Bishop's {\em
Crisis\/} text.  Their note places the original sin with {\em LEM\/}
at around 1870 (rather than Greek antiquity), when the ``flourishing
empirico-inductive tradition'' began to be replaced by the ``strictly
logico-deductive conception of pure mathematics''.  Kopell and
Stolzenberg don't hesitate to compare the {\em empirico-inductive
tradition\/} in mathematics prior to 1870, to physics, in the
following terms \cite[p.~519]{KS}:
\begin{quote}
[Mathematical] theories were theories about the phenomena, just as in
a physical theory.
\end{quote}
Similar views have been expressed by D.~Bridges~\cite{Br99}, as well
as Heyting \cite{He, He73}.  W.~Tait \cite{Ta83} argues that, unlike
intuitionism, constructive mathematics is part of classical
mathematics.  In fact, it was Frege's revolutionary logic \cite{Fre}
(see Gillies \cite{Gi92}) and other foundational developments that
created a new language and a new paradigm, transforming mathematical
foundations into fair game for further investigation, experimentation,
and criticism, including those of constructivist type.

The philosophical dilemmas in the anti-LEM sector discussed in this
section are a function of the nominalist nature of its scientific
goals.  A critique of its scientific methods appears in the next
section.

\section{A critique of the constructivist scientific method}
\label{hersh}

We would like to analyze more specifically the constructivist dilemma
with regard to the following two items:
\begin{enumerate}
\item 
the extreme value theorem, and
\item
the Hawking-Penrose singularity theorem.
\end{enumerate}

Concerning (1), note that a constructive treatment of the extreme
value theorem (EVT) by Troelstra and van Dalen \cite[p.~295]{TV}
brings to the fore the instability of the (classically unproblematic)
maximum by actually constructing a counterexample.  Such a
counterexample relies on assuming that the principle
\begin{equation*}
(a\leq 0)\vee(a\geq 0)
\end{equation*}
fails.%
\footnote{This principle is a special case of LEM; see
footnote~\ref{lpo}.}
This is a valuable insight, if viewed as a companion to classical
mathematics.%
\footnote{Some related ground on pluralism is covered in
B.~Davies~\cite{Da05}.  Mathematical phenomena such as the instability
of the extremum tend to be glossed over when approached from the
classical viewpoint; here a constructive viewpoint can provide a
welcome correction.}
If viewed as an alternative, we are forced to ponder the consequences
of the loss of the EVT.

Kronecker is sometimes thought of as the spiritual father of the
Brouwer/Bishop/Dummett tendency.  Kronecker was active at a time when
the field of mathematics was still rather compartmentalized.  Thus, he
described a 3-way partition thereof into (a) analysis, (b) geometry,
and (c) mechanics (presumably meaning mathematical physics).
Kronecker proceeded to state that it is only the analytic one-third of
mathematics that is amenable to a constructivisation in terms of the
natural numbers that ``were given to us, etc.'', but readily conceded
that such an approach is inapplicable in the remaining two-thirds,
geometry and physics.%
\footnote{\label{boniface}See Boniface and Schappacher
\cite[p.~211]{BS}.}

Nowadays mathematicians adopt a more unitary approach to the field,
and Kronecker's partition seems provincial, but in fact his caution
was vindicated by later developments, and can even be viewed as
visionary.  Consider a field such as general relativity, which in a
way is a synthesis of Kronecker's remaining two-thirds, namely,
geometry and physics.  Versions of the extreme value theorem are
routinely exploited here, in the form of the existence of solutions to
variational principles, such as geodesics, be it spacelike, timelike,
or lightlike.  At a deeper level, S.P.~Novikov \cite{No1, No2} wrote
about Hilbert's meaningful contribution to relativity theory, in the
form of discovering a Lagrangian for Einstein's equation for
spacetime.  Hilbert's deep insight was to show that general
relativity, too, can be written in Lagrangian form, which is a
satisfying conceptual insight.

A radical constructivist's reaction would be to dismiss the material
discussed in the previous paragraph as relying on LEM (needed for the
EVT), hence lacking numerical meaning, and therefore meaningless.  In
short, radical constructivism (as opposed to the liberal variety)
adopts a theory of meaning amounting to an ostrich effect%
\footnote{\label{ostrich1}Such an effect is comparable to a
traditional educator's attitude toward students' nonstandard
conceptions studied by Ely~\cite{El}, see main text in
Section~\ref{seven} around footnote~\ref{ostrich2}.}
as far as certain significant scientific insights are concerned.  A
quarter century ago, M.~Beeson already acknowledged constructivism's
problem with the calculus of variations in the following terms:
\begin{quote}
Calculus of variations is a vast and important field which lies right
on the frontier between constructive and non-constructive mathematics
\cite[p.~22]{Bee}.
\end{quote}

An even more striking example is the Hawking-Penrose singularity
theorem, whose foundational status was explored by Hellman
\cite{He98}.  The theorem relies on fixed point theorems and therefore
is also constructively unacceptable, at least in its present form.
However, the singularity theorem does provide important scientific
insight.  Roughly speaking, one of the versions of the theorem asserts
that certain natural conditions on curvature (that are arguably
satisfied experimentally in the visible universe) force the existence
of a singularity when the solution is continued backward in time,
resulting in a kind of a theoretical justification of the Big Bang.
Such an insight cannot be described as ``meaningless'' by any
reasonable standard of meaning preceding nominalist commitments.

\section{The triumvirate nominalistic reconstruction}
\label{triumvirat}

This section analyzes a nominalistic reconstruction successfully
implemented at the end of the 19th century by Cantor, Dedekind, and
Weierstrass.  The rigorisation of analysis they accomplished went
hand-in-hand with the elimination of infinitesimals; indeed, the
latter accomplishment is often viewed as a fundamental one.  We would
like to state from the outset that the main issue here is not a
nominalistic attitude on the part of our three protagonists
themselves.  Such an attitude is only clearly apparent in the case of
Cantor (see below).  Rather, we argue that the historical context in
the 1870s favored the acceptance of their reconstruction by the
mathematical community, due to a certain philosophical disposition.

Some historical background is in order.  As argued by
D.~Sherry~\cite{She87}, George Berkeley's 1734 polemical essay
\cite{Be} conflated a logical criticism and a metaphysical criticism.%
\footnote{\label{locke}Robinson distinguished between the two
criticisms in the following terms: ``The vigorous attack directed by
Berkeley against the foundations of the Calculus in the forms then
proposed is, in the first place, a brilliant exposure of their logical
inconsistencies.  But in criticizing infinitesimals of all kinds,
English or continental, Berkeley also quotes with approval a passage
in which Locke rejects the actual infinite ... It is in fact not
surprising that a philosopher in whose system perception plays the
central role, should have been unwilling to accept infinitary
entities'' \cite[p.~280-281]{Ro66}.}
In the intervening centuries, mathematicians have not distinguished
between the two criticisms sufficiently, and grew increasingly
suspicious of infinitesimals.  The metaphysical criticism stems from
the 17th century doctrine that each theoretical entity must have an
empirical counterpart/referent before such an entity can be used
meaningfully; the use of infinitesimals of course would fly in the
face of such a doctrine.%
\footnote{Namely, infinitesimals cannot be measured or perceived
(without suitable optical devices); see footnote~\ref{locke}.}

Today we no longer accept the 17th century doctrine.  However, in
addition to the metaphysical criticism, Berkeley made a poignant
logical criticism, pointing out a paradox in the definition of the
derivative.  The seeds of an approach to resolving the logical paradox
were already contained in the work of Fermat,%
\footnote{Fermat's adequality is analyzed in Section~\ref{rival1}.
Robinson modified the definition of the derivative by introducing the
standard part function, which we refer to as the Fermat-Robinson
standard part in Sections~\ref{rival1} and~\ref{rival2}.}
but it was Robinson who ironed out the remaining logical wrinkle.

Thus, mathematicians throughout the 19th century were suspicious of
infinitesimals because of a lingering influence of 17th century
doctrine, but came to reject them because of what they felt were
logical contradictions; these two aspects combined into a nominalistic
attitude that caused the triumvirate reconstruction to spread like
wildfire.

The tenor of Hobson's remarks,%
\footnote{Hobson's remarks are analyzed in footnote~\ref{hob}.}
as indeed of a majority of historians of mathematics, is that
Weierstrass's fundamental accomplishment was the elimination of
infinitesimals from foundational discourse in analysis.
Infinitesimals were replaced by arguments relying on real inequalities
and multiple-quantifier logical formulas.  

The triumvirate transformation had the effect of a steamroller
flattening a B-continuum%
\footnote{See Section~\ref{rival1}.}
into an A-continuum.  Even the ardent enthusiasts of Weierstrassian
epsilontics recognize that its practical effect on mathematical
discourse has been ``appalling"; thus, J.~Pierpont wrote as follows in
1899:
\begin{quote}
The mathematician of to-day, trained in the school of Weierstrass, is
fond of speaking of his science as `die absolut klare Wissenschaft.'
Any attempts to drag in {\em metaphysical speculations\/} are resented
with {\em indignant energy\/}.  With almost {\em painful emotions\/}
he looks back at the sorry mixture of metaphysics and mathematics
which was so common in the last century and at the beginning of this
\cite[p.~406]{Pi} [emphasis added--authors].
\end{quote}
%
%
Pierpont concludes:
\begin{quote}
The analysis of to-day is indeed a transparent science.  Built up on
the simple notion of number, its truths are the most solidly
established in the whole range of human knowledge.  It is, however,
not to be overlooked that the price paid for this clearness is {\em
appalling\/}, it is total separation from the world of our senses''
\cite[p.~406]{Pi} [emphasis added--authors].
\end{quote}
It is instructive to explore what form the ``indignant energy"
referred to by Pierpont took in practice, and what kind of rhetoric
accompanies the ``painful emotions''.  A reader attuned to 19th
century literature will not fail to recognize {\em infinitesimals\/}
as the implied target of Pierpont's epithet ``metaphysical
speculations".  

Thus, Cantor {\em published\/} a ``proof-sketch'' of a claim to the
effect that the notion of an infinitesimal is inconsistent.  By this
time, several detailed constructions of non-Archimedean systems had
appeared, notably by Stolz and du Bois-Reymond.

When Stolz published a defense of his work, arguing that technically
speaking Cantor's criticism does not apply to his system, Cantor
responded by artful innuendo aimed at undermining the credibility of
his opponents.  At no point did Cantor vouchsafe to address their
publications themselves.  In his 1890 letter to Veronese, Cantor
specifically referred to the work of Stolz and du Bois-Reymond.
Cantor refers to their work on non-Archimedean systems as not merely
an ``abomination", but a ``self contradictory and completely useless"
one.  

P.~Ehrlich \cite[p.~54]{Eh06} analyzes the errors in Cantor's
``proof'' and documents his rhetoric.

The effect on the university classroom has been pervasive.  In an
emotionally charged atmosphere, students of calculus today are warned
against taking the apparent ratio~$dy/dx$ literally.  By the time one
reaches the chain rule~$\frac{dy}{dt} = \frac{dy}{dx} \frac{dx}{dt}$,
the awkward contorsions of an obstinate denial are palpable throughout
the spectrum of the undergraduate textbooks.%
\footnote{A concrete suggestion with regard to undergraduate teaching
may be found at the end of Section~\ref{seven}.}

Who invented the real number system?  According to van der Waerden,
Simon Stevin's
\begin{quote}
general notion of a real number was accepted, tacitly or explicitly,
by all later scientists \cite[p.~69]{van}.
\end{quote}
D.~Fearnley-Sander writes that
\begin{quote}
the modern concept of real number [...] was essentially achieved by
Simon Stevin, around 1600, and was thoroughly assimilated into
mathematics in the following two centuries \cite[p.~809]{Fea}.
\end{quote}
D.~Fowler points out that
\begin{quote}
Stevin [...] was a thorough-going arithmetizer: he published, in 1585,
the first popularization of decimal fractions in the West [...]; in
1594, he desribed an algorithm for finding the decimal expansion of
the root of any polynomial, the same algorithm we find later in
Cauchy's proof of the intermediate value theorem \cite[p.~733]{Fo}.
\end{quote}
The algorithm is discussed in more detail in \cite[\S10,
p.~475-476]{Ste}.  Unlike Cauchy, who {\em halves\/} the interval at
each step, Stevin subdivides the interval into {\em ten\/} equal
parts, resulting in a gain of a new decimal digit of the solution at
every iteration of the algorithm.%
\footnote{Stevin's numbers were anticipated by E. Bonfils in 1350, see
S.~Gandz \cite{Ga}.  Bonfils says that ``the unit is divided into ten
parts which are called Primes, and each Prime is divided into ten
parts which are called Seconds, and so on into infinity''
\cite[p.~39]{Ga}.}

At variance with these historical judgments, the mathematical
community tends overwhelmingly to award the credit for constructing
the real number system to the great triumvirate,%
\footnote{See footnote~\ref{triumvirate} for the origin of this
expression.}
in appreciation of the successful extirpation of infinitesimals as a
byproduct of the Weierstrassian epsilontic formulation of analysis.

To illustrate the nature of such a reconstruction, consider Cauchy's
notion of continuity.  H.~Freudenthal notes that ``Cauchy invented our
notion of continuity'' \cite[p.~136]{Fr71a}.  Cauchy's starting point
is a description of perceptual continuity of a function in terms of
``varying by imperceptible degrees''.  Such a turn of phrase occurs
both in his letter to Coriolis of 1837, and in his 1853 text
\cite[p.~35]{Ca53}.%
\footnote{Both Cauchy's original French ``par degr\'es insensibles'',
and its correct English translation ``by imperceptible degrees'', are
etymologically related to {\em sensory perception\/}.}

Cauchy transforms perceptual continuity into a mathematical notion by
exploiting his conception of an infinitesimal as being generated by a
null sequence (see \cite{Br}).  Both in 1821 and in 1853, Cauchy
defines continuity of~$y=f(x)$ in terms of an
infinitesimal~$x$-increment resulting in an infinitesimal change
in~$y$.

The well-known nominalistic residue of the perceptual definition (a
residue that dominates our classrooms) would have~$f$ be continuous
at~$x$ if for every positive epsilon there exists a positive delta
such that if~$h$ is less than delta then~$f(x+h)-f(x)$ is less than
epsilon, namely:
\begin{equation*}
\forall \epsilon>0 \exists \delta>0 : |h| < \delta \implies
|f(x+h)-f(x)| < \epsilon.
\end{equation*}
This can hardly be said to be a hermeneutic reconstruction of Cauchy's
infinitesimal definition.  In our classrooms, are students being
dressed to perform multiple-quantifier Weierstrassian epsilontic
logical stunts, on the pretense of being taught infinitesimal
calculus?%
\footnote{One can apply here Burgess's remark to the effect that
``[t]his is educational reform in the wrong direction: away from
applications, toward entanglement in logical subtleties''
\cite[pp.~98-99]{Bu83}.}

Lord Kelvin's technician,%
\footnote{In analyzing Chihara's and Field's nominalist
reconstructions, Burgess \cite[p.~96]{Bu83} is sceptical as to the
plausibility of interpreting what Lord Kelvin's technician is saying,
in terms of tacit knowledge of such topics in foundations of
mathematics as predicative analysis and measurement theory.}
wishing to exploit the notion of continuity in a research paper, is
unlikely to be interested in 4-quantifier definitions thereof.
Regardless of the answer to such a question, the revolutionary nature
of the triumvirate reconstruction of the foundations of analysis is
evident.  If one accepts the thesis that elimination of ontological
entities called ``infinitesimals" does constitute a species of
nominalism, then the triumvirate recasting of analysis was a
nominalist project.  We will deal with Cantor and Dedekind in more
detail in Section~\ref{CD}.

\section{Cantor and Dedekind}
\label{CD}

Cantor is on record describing infinitesimals as the ``cholera
bacillus of mathematics" in a letter dated 12 december 1893, quoted in
Meschkowski \cite[p.~505]{Me} (see also Dauben \cite[p.~353]{Da95}
and~\cite[p.~124]{Da96}).  Cantor went as far as publishing a
purported ``proof" of their logical inconsistency, as discussed in
Section~\ref{triumvirat}.  Cantor may have extended numbers both in
terms of the complete ordered field of real numbers and his theory of
infinite cardinals; however, he also passionately believed that he had
not only given a logical foundation to real analysis, but also
simultaneously eliminated infinitesimals (the cholera bacillus).

Dedekind, while admitting that there is no evidence that the ``true"
continuum indeed possesses the property of completeness he championed
(see M. Moore \cite[p.~82]{Moo07}), at the same time formulated his
definition of what came to be known as Dedekind cuts, in such a way as
to rule out infinitesimals.

S.~Feferman describes Dedekind's construction of an complete ordered
line~$(\R,<)$ as follows:
\begin{quote}
Dedekind's construction of such an~$(\R, <)$ is obtained by taking it
to consist of the rational numbers together with the numbers
corresponding to all those cuts~$(X_1, X_2)$ in~$\Q$ for which~$X_1$
has no largest element and~$X_2$ has no least element, ordered in
correspondence to the ordering of cuts~$(X_1, X_2) < (Y_1, Y_2)$
when~$X_1$ is a proper subset of~$Y_1$.  Dedekind himself spoke of
this construction of~$\R$ a[s] individual cuts in~$\Q$ for which~$X_1$
has no largest element and~$X_2$ no least element as {\em the creation
of an irrational number\/},%
\footnote{Dedekind \cite{De} in 1872, translation in Ewald
\cite[p.~773]{Ew}.}
though he did not identify the numbers themselves with those cuts
\cite{Fef}.
\end{quote}
In this way the ``gappiness'' of the rationals is overcome, in
Dedekind's terminology.

Now requiring that an element of the continuum should induce a
partition of $\Q$ does not yet rule out infinitesimals.  However,
requiring that a given partition of $\Q$ should correspond to a {\em
unique\/} element of the continuum, does have the effect of ruling out
infinitesimals.  In the context of an infinitesimal-enriched
continuum, it is clear that a pair of quantities in the cluster (halo)
infinitely close to~$\pi$, for example, will define the {\em same\/}
partition of the rationals.  Therefore the clause of ``only one"
forces a collapse of the infinitesimal cluster to a single quantity,
in this case~$\pi$.

Was rigor linked to the elimination of infinitesimals?  E.~Hobson in
his retiring presidential address \cite[p.~128]{Ho} in 1902 summarized
the advances in analysis over the previous century, and went on
explicitly to make a connection between the foundational
accomplishments in analysis, on the one hand, and the elimination of
infinitesimals, on the other, by pointing out that an equivalence
class defining a real number ``is of such a character that no use is
made in it of infinitesimals,''%
\footnote{See footnote~\ref{hob} for more details on Hobson.}
suggesting that Hobson viewed them as logically inconsistent (perhaps
following Cantor or Berkeley).  The matter of rigor will be analyzed
in more detail in the next section.

\section{A critique of the Weierstrassian scientific method}
\label{six}

In criticizing the nominalistic aspect of the Weierstrassian
elimination of infinitesimals, does one neglect the mathematical
reasons why this was considered desirable?

The stated goal of the triumvirate program was mathematical rigor.
Let us examine the meaning of mathematical rigor.  Conceivably rigor
could be interpreted in at least the following four ways, not all of
which we endorse:

\begin{enumerate}
\item
it is a shibboleth that identifies the speaker as belonging to a clan
of professional mathematicians;
\item
it represents the idea that as the field develops, its practitioners
attain greater and more conceptual understanding of key issues, and
are less prone to error;
\item
it represents the idea that a search for greater correctness in
analysis inevitably led Weierstrass to epsilontics in the 1870s;
\item
it refers to the establishment of ultimate foundations for mathematics
by Cantor, eventually explicitly expressed in axiomatic form by
Zermelo and Fraenkel.%
\footnote{\label{ist}It would be interesting to investigate the role
of the Zermelo--Frankel axiomatisation of set theory in cementing the
nominalistic disposition we are analyzing.  Keisler points out that
``the second- and higher-order theories of the real line depend on the
underlying universe of set theory [...] Thus the {\em properties\/} of
the real line are not uniquely determined by the axioms of set
theory'' \cite[p.~228]{Ke94} [emphasis in the original--the authors].
He adds: ``A set theory which was not strong enough to prove the
unique existence of the real line would not have gained acceptance as
a mathematical foundation'' \cite[p.~228]{Ke94}.  Edward Nelson
\cite{Ne} has developed an alternative axiomatisation more congenial
to infinitesimals.}
\end{enumerate}

Item (1) may be pursued by a fashionable academic in the social
sciences, but does not get to the bottom of the issue.  Meanwhile,
item (2) could apply to any exact science, and does not involve a
commitment as to which route the development of mathematics may have
taken.  Item (2) could be supported by scientists in and outside of
mathematics alike, as it does not entail a commitment to a specific
destination or ultimate goal of scientific devepment as being
pre-determined and intrinsically inevitable.

On the other hand, the actual position of a majority of professional
mathematicians today corresponds to items (3) and (4).  The crucial
element present in (3) and (4) and absent in (2) is the postulation of
a specific outcome, believed to be the inevitable result of the
development of the discipline.  Challenging such a {\em belief\/}
appears to be a radical proposition in the eyes of a typical
professional mathematician, but not in the eyes of scientists in
related fields of the exact sciences.  It is therefore particularly
puzzling that (3) and (4) should be accepted without challenge by a
majority of historians of mathematics, who tend to toe the line on the
mathematicians' {\em belief\/}.  It is therefore necessary to examine
such a belief, which, as we argue, stems from a particular
philosophical disposition akin to nominalism.

Could mathematical analysis have followed a different path of
development?  In an intriguing text published a decade ago,
Pourciau~\cite{Po} examines the foundational crisis of the 1920s and
the Brouwer--Hilbert controversy, and argues that Brouwer's view may
have prevailed had Brouwer been more of an$\ldots$ Errett Bishop.%
\footnote{More specifically, Pourciau would have wanted Brouwer to
stop ``wasting time'' on free choice sequences and the continuum, and
to focus instead on developing analysis on a constructive footing
based on~$\N$.}
While we are sceptical as to Pourciau's main conclusions, the
unmistakable facts are as follows:
\begin{enumerate}
\item
a real struggle did take place;
\item
some of the most brilliant minds at the time did side with Brouwer, at
least for a period of time (e.g., Hermann Weyl);
\item
the battle was won by Hilbert not by mathematical means alone but also
political ones, such as maneuvering Brouwer out of a key editorial
board;
\item
while retroactively one can offer numerous reasons why Hilbert's
victory might have been inevitable, this was not at all obvious at the
time.
\end{enumerate}

We now leap back a century, and consider a key transitional figure,
namely Cauchy.  In 1821, Cauchy defined continuity of~$y=f(x)$ in
terms of ``an infinitesimal~$x$-increment corresponding to an
infinitesimal~$y$-increment''.%
\footnote{See Section~\ref{battle} for more details on Cauchy's
definition.}
Many a practicing mathematician, brought up on an alleged
``Cauchy-Weierstrass~$\epsilon, \delta$'' tale, will be startled by
such a revelation.  The textbooks and the history books routinely
obfuscate the nature of Cauchy's definition of continuity.

Fifty years before Weierstrass, Cauchy performed a hypostatisation by
encapsulating a variable quantity tending to zero, into an
individual/atomic entity called ``an infinitesimal''.


Was the naive traditional ``definition" of the infinitesimal blatantly
self-contradictory?  We argue that it was not.  Cauchy's definition in
terms of null sequences is a reasonable definition, and one that
connects well with the sequential approach of the ultrapower
construction.%
\footnote{See discussion around formula~\eqref{infinitesimal} in
Section~\ref{rival2} for more details.}
Mathematicians viewed infinitesimals with deep suspicion due in part
to a conflation of two separate criticisms, the logical one and the
metaphysical one, by Berkeley, see Sherry~\cite{She87}.  Thus, the
emphasis on the elimination of infinitesimals in the traditional
account of the history of analysis is misplaced.  

Could analysis could have developed on the basis of infinitesimals?%
\footnote{The issue of alternative axiomatisations is discussed in
footnote~\ref{ist}.}
Continuity, as all fundamental notions of analysis, can be defined,
and were defined, by Cauchy in terms of infinitesimals.  Epsilontics
could have played a secondary role of clarifying whatever technical
situations were too awkward to handle otherwise, but arguably they
needn't have {\em replaced\/} infinitesimals.  As far as the issue of
rigor is concerned, it needs to be recognized that Gauss and Dirichlet
published virtually error-free mathematics before Weierstrassian
epsilontics, while Weierstrass himself was not protected by
epsilontics from publishing an erroneous paper by S. Kovalevskaya (the
error was found by Volterra, see \cite[p.~568]{Ra}).  As a scientific
discipline develops, its practitioners gain a better understanding of
the conceptual issues, which helps them avoid errors.  But assigning a
singular, oracular, and benevolent role in this to epsilontics is
philosophically naive.  The proclivity to place the blame for errors
on infinitesimals betrays a nominalistic disposition aimed agaist the
{\em ghosts of departed quantities\/}, already dubbed ``{\em
charlatanerie\/}'' by d'Alembert \cite{dal} in 1754.

\section{A question session}
\label{seven}

The following seven questions were formulated by R.~Hersh, who also
motivated the author to present the material of Section~\ref{hersh},
as well as that of Section~\ref{six}.

\begin{question}
Was a nominalistic viewpoint motivating the triumvirate project?
\end{question}

Answer.  We argue that the answer is affirmative, and cite two items
as evidence:
\begin{enumerate}
\item
Dedekind's cuts and the ``essence of continuity", and
\item
Cantor's tooth-and-nail fight against infinitesimals.
\end{enumerate}

Concerning (1), mathematicians widely believed that Dedekind
discovered such an essence.  What is meant by the essence of
continuity in this context is the idea that a pair of ``cuts" on the
rationals are identical if and only if the pair of numbers defining
them are equal.  Now the ``if" part is unobjectionable, but the almost
reflexive ``only if" part following it has the effect of a steamroller
flattening the B-continuum%
\footnote{See Section~\ref{rival1}.}
into an A-continuum.  Namely, it collapses each monad (halo, cluster)
to a point, since a pair of infinitely close (adequal) points
necessarily define the same cut on the rationals (see Section~\ref{CD}
for more details).  The fact that the steamroller effect was gladly
accepted as a near-axiom is a reflection of a nominalistic attitude.

Concerning (2), Cantor not only published a ``proof-sketch" of the
non-existence of infinitesimals, he is on record calling them an
``abomination" as well as the ``cholera bacillus" of mathematics.
When Stolz meekly objected that Cantor's ``proof" does not apply to
his system, Cantor responded by the ``abomination" remark (see
Section~\ref{CD} for more details).  Now Cantor's proof contains an
error that was exhastively analyzed by Ehrlich \cite{Eh06}.  As it
stands, it would ``prove'' the non-existence of the surreals!
Incidentally, Ehrlich recently proved that ``maximal" surreals are
isomorphic to ``maximal" hyperreals.  Can Cantor's attitude be
considered as a philosophical predisposition to the detriment of
infinitesimals?

\begin{question}
It has been written that Cauchy's concern with clarifying the
foundations of calculus was motivated by the need to teach it to
French military cadets.
\end{question}

Answer.  Cauchy did have some tensions with the management of the {\em
Ecole Polytechnique\/} over the teaching of infinitesimals between
1814 and 1820.  In around 1820 he started using them in earnest in
both his textbooks and his research papers, and continued using them
throughout his life, well past his teaching stint at the {\em
Ecole\/}.  Thus, in his 1853 text \cite{Ca53} he reaffirms the
infinitesimal definition of continuity he gave in his 1821
textbook~\cite{Ca21}.

\begin{question}
Doesn't reasoning by infinitesimals require a deep intuition that is
beyond the reach of most students?
\end{question}

Kathleen Sullivan's study \cite{Su} from 1976 shows that students
enrolled in sections based on Keisler's textbook end up having a
better conceptual grasp of the notions of calculus than control groups
following the standard approach.  Two years ago, I taught
infinitesimal calculus to a group of 25 freshmen.  I also had to train
the TA who was new to the material.  According to persistent reports
from the TA, the students have never been so excited about learning
calculus.  On the contrary, it is the multiple-quantifier
Weierstrassian epsilontic logical stunts that our students are dressed
to perform (on pretense of being taught infinitesimal calculus) that
are beyond their reach.  In an ironic commentary on the nominalistic
ethos reigning in our departments, not only was I relieved of my
teaching this course the following year, but the course number itself
was eliminated.

\begin{question}
It may true that ``epsilontics" is in practice repugnant to many
students.  But the question is whether an issue that is really a
matter of technical mathematics related to pedagogy is being
misleadingly presented as a question of high metaphysics.
\end{question}

Answer.  Berkeley turned this into a metaphysical debate.  Generations
of mathematicians have grown up thinking of it as a metaphysical
debate.  Such a characterisation is precisely what we contest.

\begin{question}
Isn't ``a positive number smaller than all positive numbers"
self-contradictory?  Is a phrase such as ``I am smaller than myself",
intelligible?
\end{question}

Answer.  Both Carnot and Cauchy say that an infinitesimal is generated
by a variable quantity that becomes smaller than any fixed quantity.
No contradiction here.  The otherwise excellent study by
Ehrlich~\cite{Eh06} contains a curious slip with regard to Poisson.
Poisson describes infinitesimals as being ``less than any given
magnitude of the same nature'' \cite[p.~13-14]{Poi} (the quote is
reproduced in Boyer \cite[p.~283]{Boy}).  Ehrlich inexplicably omits
the crucial modifier ``given'' when quoting Poisson in footnote~133 on
page~76 of \cite{Eh06}.  Based on the incomplete quote, Ehrlich
proceeds to agree with Veronese's assessment (of Poisson) that
``[t]his proposition evidently contains a contradiction in terms"
\cite[p.~622]{Ve91}.  Our assessment is that Poisson's definition is
in fact perfectly consistent.

\begin{question}
Infinitesimals were one thorny issue.  Didn't it take the modern
theory of formal languages to untangle that?
\end{question}

Answer.  Not exactly.  A long tradition of technical work in
non-Archimedean continua starts with Stolz and du Bois-Reymond,
Levi-Civita, Hilbert, and Borel, see Ehrlich~\cite{Eh06}.  The
tradition continues uninterruptedly until Hewitt constructs the
hyperreals in 1948.  Then came {\L}o{\'s}'s theorem whose consequence
is a transfer principle, which is a mathematical implementation of the
heuristic ``law of continuity" of Leibniz (``what's true in the finite
domain should remain true in the infinite domain").  What {\L}o{\'s}
and Robinson untangled was the transfer principle.  Non-Archimedean
systems had a long history prior to these developments.

\begin{question}
There is still a pedagogical issue.  I do understand that Keisler's
calculus book is teachable.  But this says nothing about the
difficulty of teaching calculus in terms of infinitesimals back around
1800.  Keisler has Robinson's non-standard analysis available, as a
way to make sense of infinitesimals.  Cauchy did not.  Do you believe
that Cauchy used a definition of infinitesimal in terms of a null
sequence of rationals (or reals) in teaching introductory calculus?
\end{question}

Answer.  The historical issue about Cauchy is an interesting one.
Most of his course notes from the period 1814-1820 have been lost.
His predecessor at the {\em Ecole Polytechnique\/}, L.~Carnot, defined
infinitesimals exactly the same way as Cauchy did, but somehow is
typically viewed by historians as belonging to the old school as far
as infinitesimals are concerned (and criticized for his own version of
the ``cancellation of errors'' argument originating with Berkeley).
As far as Cauchy's textbooks from 1821 onward indicate, he declares at
the outset that infinitesimals are an indispensable foundational tool,
defines them in terms of null sequences (more specifically, ``a
variable quantity becomes an infinitesimal"), defines continuity in
terms of infinitesimals, defines his ``Dirac" delta function in terms
of infinitesimals (see \cite{Lau92}), defines infinitesimals of
arbitrary {\em real\/} order in \cite[p. 281]{Ca29}, anticipating
later work by Stolz, du Bois-Reymond, and others.

The following eight questions were posed by Martin Davis.

\begin{question}
How would you answer the query: how do you define a ``null sequence"?
Aren't you back to epsilons?
\end{question}

Answer.  Not necessarily.  One could define it, for example, in terms
of ``only finitely many terms outside each given separation from
zero''.  While epsilontics has important applications, codifying the
notion of a null sequence is not one of them.  Epsilontics is helpful
when it comes to characterizing a Cauchy sequence, {\em if one does
not yet know the limiting value\/}.  If one does know the limiting
value, as in the case of a null sequence, a multiple-quantifier
epsilontic formulation is no clearer than saying that all the terms
eventually get arbitrarily small.  To be more specific, if one
describes a Cauchy sequence by saying that ``terms eventually get
arbitrarily close to each other'', the ambiguity can lead and has led
to errors, though not in Cauchy (the sequences are rightfully named
after him as he was aware of the trap).  Such an ambiguity is just not
there as far as null sequences are concerned.  Giving an epsilontic
definition of a null sequence does not increase understanding and does
not decrease the likelihood of error.  A null sequence is arguably a
notion that's no more complex than multiple-quantifier epsilontics,
just as the natural numbers are no more complex than the set-theoretic
definition thereof in terms of~$0=\emptyset, \; 1=\{\emptyset \}, \; 2
= \{\emptyset, \{ \emptyset \} \}, \; \ldots$ which requires
infinitely many set-theoretic types to make the point.  Isn't a
natural number a more primitive notion?

\begin{question}
You evoke Cauchy's use of variable quantities.  But whatever is a
``variable quantity"?
\end{question}

The concept of variable quantity was not clearly defined by
mathematicians from Leibniz and l'Hopital onwards, and is today
considered a historical curiosity.  Cauchy himself sometimes seems to
think they take discrete values (in his 1821 text \cite{Ca21}) and
sometimes continuous (in his 1823 text \cite{Ca23}).  Many historians
agree that in 1821 they were discrete sequences, and Cauchy himself
gives explicit examples of sequences.  Now it is instructive to
compare such variable quantities to the procedures advocated by the
triumvirate.  In fact, the approach can be compared to Cantor's
construction of the real numbers.  A real number is a Cauchy sequence
of rational numbers, modulo an equivalence relation.  A sequence,
which is not an individual/atomic entity, comes to be viewed as an
atomic entity by a process which in triumvirate lexicon is called an
equivalence relation, but in much older philosophical terminology is
called {\em hypostatisation\/}.  In education circles, researchers
tend to use terms such as {\em encapsulation\/}, {\em procept\/}, and
{\em reification\/}, instead.  As you know, the ultrapower
construction is a way of hypostatizing a hyperreal out of sequences of
reals.  As far as Cauchy's competing views of a variable quantity as
discrete (in 1821) or continuous (in 1823), they lead to distinct
implementations of a B-continuum in Hewitt (1948), when a
``continuous'' version of the ultrapower construction was used, and in
Luxemburg (1962), where the discrete version was used (a more recent
account of the latter is in Goldblatt~\cite{Go}).

\begin{question}
Isn't the notion of a variable quantity a pernicious notion that makes
time an essential part of mathematics?
\end{question}

Answer.  I think you take Zeno's paradoxes too seriously.  I
personally don't think there is anything wrong with involving time in
mathematics.  It has not led to any errors as far as I know, {\em
pace\/} Zeno.

\begin{question}
Given what relativity has taught us about time, is it a good idea to
involve time in mathematics?
\end{question}

Answer.  What did relativity teach us about time that would make us
take time out of mathematics?  That time is relative?  But you may be
confusing metaphysics with mathematics.  Time that's being used in
mathematics is not an exact replica of physical time.  We may still be
influenced by 17th century doctrine according to which every
theoretical entity must have an empirical counterpart/referent.  This
is why Berkeley was objecting to infinitesimals (his metaphysical
criticism anyway).  I would put time back in mathematics the same way
I would put infinitesimals back in mathematics.  Neither concept is
under any obligation of corresponding to an empirical referent.

\begin{question}
Didn't the triumvirate show us how to prove the existence of a
complete ordered field?
\end{question}

Answer.  Simon Stevin had already made major strides in defining the
real numbers, represented by decimals.  Some essential work needed to
be done, such as the fact that the usual operations are well defined.
This was done by Dedekind, see Fowler \cite{Fo}.  But Stevin numbers
themselves were several centuries older, even though they go under the
soothing name of numbers so {\em real\/}.%
\footnote{This theme is developed in more detail in
Section~\ref{triumvirat}.}

\begin{question}
My NSA book \cite{Dav77} does it by forming the quotient of the ring
of finite hyper-rational numbers by the ideal of
infinitesimals$\ldots$
\end{question}

The remarkable fact is that this construction is already anticipated
by K\"astner
%
%
(a contemporary of Euler's) in the following terms: ``If one
partitions~$1$ without end into smaller and smaller parts, and takes
larger and larger collections of such little parts, one gets closer
and closer to the irrational number without ever attaining it''.
K\"astner concludes:
\begin{quote}
``Therefore one can view it as an infinite collection of infinitely
small parts'' \cite{Ka}, cited by Cousquer \cite{Cou}.%
\footnote{K\"astner's suggestion is implemented by the surjective
leftmost vertical arrow in our Figure~\ref{helpful} in
Section~\ref{rival2}.}
\end{quote}

\begin{question}
$\ldots$ but that construction was not available to the earlier
generations.
\end{question}

But Stevin numbers were.  They kept on teaching analysis in France
throughout the 1870s without any need for ``constructing" something
that had already been around for a century before Leibniz, see
discussion in Laugwitz \cite[p.~274]{Lau00}.

\begin{question}
Aren't you conflating the problem of rigorous foundation with how to
teach calculus to beginners?
\end{question}

As far as rigorous foundations are concerned, alternative foundations
to ZF have been developed that are more congenial to infinitesimals,
such as Edward Nelson's~\cite{Ne}.  Mathematicians are accustomed to
thinking of ZF as ``the foundations".  It needs to be recognized that
this is a philosophical assumption.  The assumption can be a
reflection of a nominalist mindframe.

The example of a useful application of infinitesimals analyzed at the
end of this section is quite elementary.  A more advanced example is
the elegant construction of the Haar measure in terms of counting
infinitesimal neighborhoods, see Goldblatt \cite{Go}.  Even more
advanced examples such as the proof of the invariant subspace
conjecture are explained in your book~\cite{Dav77}.  For an
application to the Bolzmann equation, see Arkeryd~\cite{Ar81, Ar05}.

\begin{figure}
\includegraphics[height=2in]{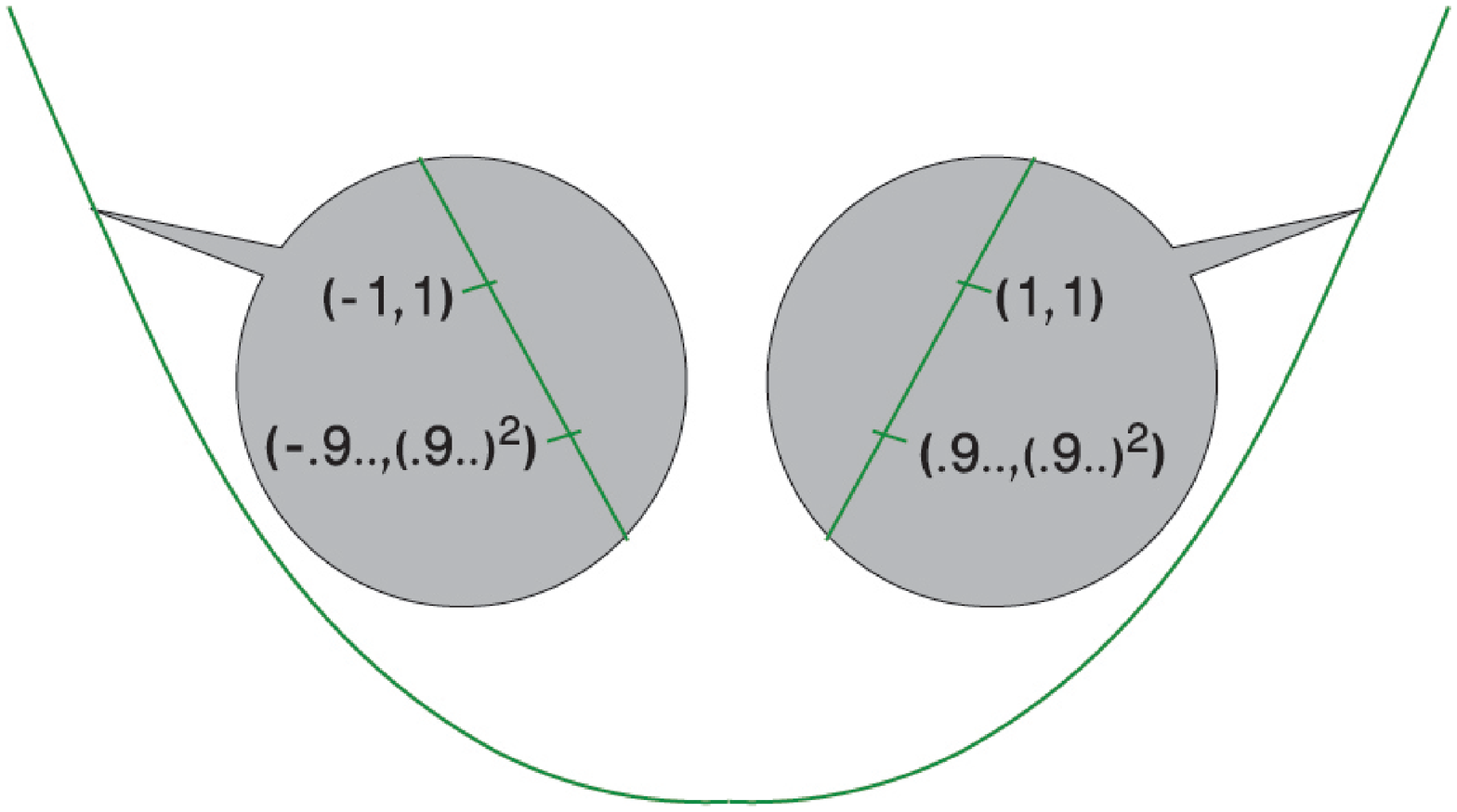}
\caption{\textsf{Differentiating~$y=f(x)=x^2$ at~$x=1$ yields
$\tfrac{\Delta y}{\Delta x} = \tfrac{f(.9..) - f(1)}{.9..-1} =
\tfrac{(.9..)^2 - 1}{.9..-1} = \tfrac{(.9.. - 1)(.9.. + 1)}{.9..-1} =
.9.. + 1 \approx 2.$ Here~$\approx$ is the relation of being
infinitely close (adequal).  Hyperreals of the form~$.9..$ are
discussed in \cite{KK2}}}
\label{jul10}
\end{figure}

As a concrete example of what consequences a correction of the
nominalistic triumvirate attitude would entail in the teaching of the
calculus, consider the problem of the {\em unital evaluation\/} of the
decimal~$.999\ldots$, i.e., its evaluation to the {\em unit\/}
value~$1$.  Students are known overwhelmingly to believe that the
number~$a=.999\dots$ falls short of~$1$ by an infinitesimal amount.  A
typical instructor believes such student intuitions to be erroneous,
and seeks to inculcate the unital evaluation of~$a$.  An alternative
approach was proposed by Ely \cite{El} and Katz \& Katz~\cite{KK1}.
Instead of refuting student intuitions, an instructor could build upon
them to calculate the derivative of~$y=x^2$ at~$x=1$ by choosing an
infinitesimal~$\Delta x=a-1$ and showing that
\begin{equation*}
\frac{\Delta y}{\Delta
x}=\frac{a^2-1^2}{a-1}=\frac{(a-1)(a+1)}{a-1}=a+1
\end{equation*}
is infinitely close (adequal) to~$2$, yielding the desired value
without either epsilontics, estimates, or limits, see
Figure~\ref{jul10}.  Here~$a$ is interpreted as an extended decimal
string with an infinite hypernatural's worth of~$9$s, see~\cite{Li}.
Instead of building upon student intuition, a typical calculus course
seeks to flatten it into the ground by steamrolling the B-continuum
into the A-continuum, see Katz and Katz \cite{KK1, KK2}.  A nominalist
view of what constitutes an allowable number system has produced an
ostrich effect%
\footnote{\label{ostrich2}Such an effect is comparable to a
constructivist's reaction to the challenge of meaningful applications
of a post-LEM variety, see main text in Section~\ref{hersh} around
footnote~\ref{ostrich1}.}
whereby mathematics educators around the globe have failed to
recognize the legitimacy, and potency, of students' nonstandard
conceptions of~$.999\ldots$, see Ely~\cite{El} for details.

\section{The battle for Cauchy's lineage}
\label{battle}

This section analyzes the reconstruction of Cauchy's foundational work
in analysis usually associated with J.~Grabiner, and has its sources
in the work of C.~Boyer.  A critical analysis of the traditional
approach may be found in Hourya Benis Sinaceur's article~\cite{Si}
from 1973.  To place such work in a historical perspective, a minimal
chronology of commentators on Cauchy's foundational work in analysis
would have to mention F.~Klein's observation in 1908 that
\begin{quote}
since Cauchy's time, the words {\em infinitely small\/} are used in
modern textbooks in a somewhat changed sense.  One never says, namely,
that a quantity {\em is\/} infinitely small, but rather that it {\em
becomes\/} infinitely small'' \cite[p.~219]{Kl}.
\end{quote}
Indeed, Cauchy's starting point in defining an infinitesimal is a null
sequence (i.e., sequence tending to zero), and he repeatedly refers to
such a null sequence as {\em becoming\/} an infinitesimal.

P.~Jourdain's detailed 1913 study \cite{Jo} of Cauchy is characterized
by a total {\em absence\/} of any claim to the effect that Cauchy may
have based his notion of infinitesimal, on limits.

C.~Boyer quotes Cauchy's definition of continuity as follows: 
\begin{quote}
the function~$f$ is continuous within given limits if between these
limits an infinitely small increment~$i$ in the variable~$x$ produces
always an infinitely small increment,~$f(x+i)-f(x)$, in the function
itself \cite[p.~277]{Boy}.
\end{quote}
Next, Boyer proceeds to {\em interpret\/} Cauchy's definition of
continuity as follows: ``The expressions infinitely small {\em are
here to be understood\/} [...]  in terms of [...]  limits: i.e.,
$f(x)$ is continuous within an interval if the limit of the variable
$f(x)$ as~$x$ approaches~$a$ is~$f(a)$, for any value of~$a$ within
this interval'' [emphasis added--authors].  Boyer feels that
infinitesimals {\em are to be understood\/} in terms of limits.  Or
perhaps they are to be understood otherwise?

In 1967, A.~Robinson discussed the place of infinitesimals in Cauchy's
work.  He pointed out that ``the assumption that [infinitesimals]
satisfy the same laws as the ordinary numbers, which was stated
explicitly by Leibniz, was rejected by Cauchy as unwarranted''.  Yet,
\begin{quote}
Cauchy's professed opinions in these matters notwithstanding, he did
in fact treat infinitesimals habitually as if they were ordinary
numbers and satisfied the familiar rules of arithmetic%
\footnote{Freudenthal, similarly, notes a general tendency on Cauchy's
part not to play by the rules: ``Cauchy was rather more flexible than
dogmatic, for more often than not he sinned against his own precepts''
\cite[p.~137]{Fr71a}.  Cauchy's irrevent attitude extended into the
civic domain, as Freudenthal reports the anecdote dealing with
Cauchy's stint as social worker in the town of Sceaux: ``he spent his
entire salary for the poor of that town, about which behavior he
reassured the mayor: `Do not worry, it is only my salary; it is not my
money, it is the emperor's' '' \cite[p.~133]{Fr71a}.}
\cite[p.~36]{Ro67}, \cite[p.~545]{Ro79}.
\end{quote}

T.~Koetsier remarked that, had Cauchy wished to extend the domain of
his functions to include infinitesimals, he would no doubt have
mentioned how exactly the functions are to be so extended.%
\footnote{Here Koetsier asks: ``If Cauchy had really wanted to
consider functions defined on sets of infinitesimals, isn't it then
highly improbable that he would not have explicitly said so?''
\cite[p.~90]{Ko}.}
Beyond the observation that Cauchy did, in fact, make it clear that
such an extension is to be carried out term-by-term,%
\footnote{Namely, an infinitesimal being generated by a null sequence,
we evaluate~$f$ at it by applying~$f$ to each term in the sequence.
Br\aa ting \cite{Br} analyzes Cauchy's use of the particular
sequence~$x=\frac{1}{n}$ in \cite{Ca53}.}
Koetsier's question prompts a similar query: had Cauchy wished to base
his calculus on limits, he would no doubt have mentioned something
about such a foundational stance.  Instead, Cauchy emphasized that in
founding analysis he was unable to avoid elaborating the fundamental
properties of {\em infinitely small quantities\/}, see \cite{Ca21}.
No mention of a foundational role of limits is anywhere to be found in
Cauchy, unlike his would-be modern interpreters.

L.~Sad {\em et al\/} have pursued this matter in detail in \cite{STB},
arguing that what Cauchy had in mind was a prototype of an ultrapower
construction, where the equivalence class of a null sequence indeed
produces an infinitesimal, in a suitable set-theoretic framework.%
\footnote{See formula~\eqref{infinitesimal} in Section~\ref{rival2}.}

To summarize, a post-Jourdain nominalist reconstruction of Cauchy's
infinitesimals, originating no later than Boyer, reduces them to a
Weierstrassian notion of limit.  To use Burgess' terminology borrowed
from linguistics, the Boyer-Grabiner interpretation
\begin{quote}
becomes the hypothesis that certain {\em noun phrases\/} [in the
present case, infinitesimals] in the surface structure are without
counterpart in the deep structure \cite[p.~97]{Bu83}.
\end{quote}
Meanwhile, a rival school of thought places Cauchy's continuum firmly
in the line of infinitesimal-enriched continua.  The ongoing debate
between rival visions of Cauchy's continuum echoes Felix Klein's
sentiment reproduced above.%
\footnote{\label{cauchy}See discussion of Klein in the main text
around footnote~\ref{klein}.  The two rival views of Cauchy's
infinitesimals have been pursued by historians, mathematicians, and
philosophers, alike.  The bibliography in the subject is vast.  The
most detailed statement of Boyer's position may be found in Grabiner
\cite{Gr81}.  Robinson's perspective was developed most successfully
by D.~Laugwitz ~\cite{Lau89} in 1989, and by K.~Br\aa ting~\cite{Br}
in 2007.}

Viewed through the lens of the dichotomy introduced by Burgess, it
appears that the traditional Boyer-Grabiner view is best described as
a hermeneutic, rather than revolutionary, nominalistic reconstruction
of Cauchy's foundational work.

Cauchy's definition of continuity in terms of infinitesimals has been
a source of an on-going controversy, which provides insight into the
nominalist nature of the Boyer-Grabiner reconstruction.  Many
historians have interpreted Cauchy's definition as a
proto-Weierstrassian definition of continuity in terms of limits.
Thus, Smithies \cite[p.~53, footnote~20]{Sm} cites the {\em page\/} in
Cauchy's book where Cauchy gave the infinitesimal definition, but goes
on to claim that the concept of {\em limit\/} was Cauchy's ``essential
basis'' for his concept of continuity \cite[p.~58]{Sm}.  Smithies
looked in Cauchy, saw the infinitesimal definition, and went on to
write in his paper that he saw a limit definition.  Such automated
translation has been prevalent at least since Boyer
\cite[p.~277]{Boy}.  Smithies cites chapter and verse in Cauchy where
the latter gives an infinitesimal definition of continuity, and
proceeds to claim that Cauchy gave a modern one.  Such awkward
contortions are a trademark of a nominalist.  In the next section, we
will examine the methodology of nominalistic Cauchy scholarship.

\section{A subtle and difficult idea of adequality}

The view of the history of analysis from the 1670s to the 1870s as a
2-century triumphant march toward the yawning heights of the rigor%
\footnote{or {\em rigor mathematicae\/}, as D.~Sherry put it in
\cite{She87}}
of Weierstrassian epsilontics has permeated the very language
mathematicians speak today, making an alternative account nearly
unthinkable.  A majority of historians have followed suit, though some
truly original thinkers differed.  These include C.~S.~Peirce, Felix
Klein, N.~N.~Luzin~\cite{Luz31}, Hans Freudenthal, Robinson,
Lakatos~\cite{La}, Laugwitz, Teixeira \cite{STB}, and Br\aa ting
\cite{Br}.

Meanwhile, J.~Grabiner offered the following reflection on the subject
of George Berkeley's criticism of infinitesimal calculus:
\begin{quote}
[s]ince an adequate response to Berkeley's objections would have
involved recognizing that an equation involving limits is a shorthand
expression for a sequence of inequalities---a subtle and difficult
idea---no eighteenth century analyst gave a fully adequate answer to
Berkeley \cite[p.~189]{Grab}.
\end{quote}
This is an astonishing claim, which amounts to reading back into
history, feedback-style, developments that came much later.%
\footnote{\label{GG}Grattan-Guinness enunciates a historical
reconstruction project in the name of H.~Freudenthal~\cite{Fr71b} in
the following terms:
\begin{quote}
``it is mere feedback-style ahistory to read Cauchy (and
contemporaries such as Bernard Bolzano) as if they had read
Weierstrass already.  On the contrary, their own pre-Weierstrassian
muddles need historical reconstruction \cite[p.~176]{Grat04}.
\end{quote}
The term ``muddle'' refers to an irreducible ambiguity of historical
mathematics such as Cauchy's sum theorem of 1821.}
Such a claim amounts to postulating the inevitability of a triumphant
march, from Berkeley onward, toward the radiant future of
Weierstrassian epsilontics (``sequence of inequalities---a subtle and
difficult idea'').  The claim of such inevitability in our opinion is
an assumption that requires further argument.  Berkeley was, after
all, attacking the coherence of {\em infinitesimals\/}.  He was not
attacking the coherence of some kind of incipient form of
Weierstrassian epsilontics and its inequalities.  Isn't there a
simpler answer to Berkeley's query, in terms of a distinction between
``variable quantity'' and ``given quantity'' already present in
l'H\^opital's textbook at the end of the 17th century?  The missing
ingredient was a way of relating a variable quantity to a given
quantity, but that, too, was anticipated by Pierre de Fermat's concept
of adequality, as discussed in Section~\ref{rival1}.

We will analyze the problem in more detail from the 19th century,
pre-Weierstrass, viewpoint of Cauchy's textbooks.  In Cauchy's world,
a variable quantity~$q$ can have a ``limiting'' fixed quantity~$k$,
such that the difference~$q-k$ is infinitesimal.  Consider Cauchy's
decomposition of an arbitrary infinitesimal of order~$n$ as a sum
\begin{equation*}
k\alpha^n(1+\epsilon)
\end{equation*}
(see \cite[p.~28]{Ca21}), where~$k$ is fixed nonzero, whereas
$\epsilon$ is a variable quantity representing an infinitesimal.  If
one were to set~$n=0$ in this formula, one would obtain a
representation of an an arbitrary finite quantity~$q$, as a sum
\begin{equation*}
q=k+k\epsilon.
\end{equation*}
If we were to suppress the infinitesimal part~$k\epsilon$, we would
obtain ``the standard part''~$k$ of the original variable
quantity~$q$.  In the terminology of Section~\ref{rival1} we are
dealing with a passage from a finite point of a B-continuum, to the
infinitely close (adequal) point of the A-continuum, namely passing
from a variable quantity to its limiting constant (fixed, given)
quantity.

Cauchy had the means at his disposal to resolve Berkeley's query, so
as to solve the logical puzzle of the definition of the derivative in
the context of a B-continuum.  While he did not resolve it, he did not
need the subtle and difficult idea of Weierstrassian epsilontics;
suggesting otherwise amounts to feedback-style ahistory.

This reader was shocked to discover, upon his first reading of
chapter~6 in Schubring \cite{Sch}, that {\em Schubring is not aware of
the fact that Robinson's non-standard numbers are an extension of the
real numbers\/}.

Consider the following three consecutive sentences from Schubring's
chapter~6: 
\begin{quote}
``[A]~[Giusti's 1984 paper] spurred Laugwitz to even more detailed
attempts to banish the error and confirm that Cauchy had used
hyper-real numbers.

\noindent
[B] On this basis, he claims, the errors vanish and the theorems
become correct, or, rather, they always were correct (see Laugwitz
1990, 21).

\noindent
[C] In contrast to Robinson and his followers, Laugwitz (1987) assumes
that Cauchy did not use nonstandard numbers in the sense of NSA, but
that his {\em infiniment petits\/} were infinitesimals representing an
extension of the field of real numbers'' \cite[p.~432]{Sch}.
\end{quote}
These three sentences, which we have labeled [A], [B], and~[C], tell a
remarkable story that will allow us to gauge Schubring's exact
relationship to the subject of his speculations.  What interests us
are the first sentence [A] and the last sentence~[C].  Their literal
reading yields the following four views put forth by Schubring:
(1)~Laugwitz intepreted Cauchy as using hyperreal numbers (from
sentence~[A]); (2)~Robinson assumed that Cauchy used ``nonstandard
numbers in the sense of NSA" (from sentence [C]); (3)~Laugwitz
disagreed with Robinson on the latter point (from sentence~[C]); (4)
Laugwitz interpreted Cauchy as using an extension of the field of real
numbers (from sentence~[C]).  Taken at face value, items~(1) and (4)
together would logically indicate that~(5) Laugwitz interpreted Cauchy
as using the hyperreal extension of the reals; moreover, if, as
indicated in item~(3), Laugwitz disagreed with Robinson, then it would
logically follow that (6)~Robinson interpreted Cauchy as {\em not\/}
using the hyperreal extension of the reals; as to the question what
number system Robinson {\em did\/} attribute to Cauchy, item~(2) would
indicate that~(7) Robinson used, not Laugwitz's hyperreals, but rather
``nonstandard numbers in the sense of NSA".

We hasten to clarify that all of the items listed above are
incoherent.  Indeed, Robinson's ``non-standard numbers" and the
hyperreals are one and the same number system (see
Section~\ref{rival2} for more details; Robinson's approach is actually
more general than Hewitt's hyperreal fields).  Meanwhile, Laugwitz's
preferred system is a different system altogether, called
Omega-calculus.  We gather that Schubring literally does not know what
he is writing about when he takes on Robinson and Laugwitz.

A reader interested in an introduction to Popper and fallibilism need
look no further than chapter~6 of Schubring \cite{Sch}, who comments
on
\begin{quote}
the enthusiasm for revising traditional beliefs in the history of
science and reinterpreting the discipline from a theoretical,
epistemological perspective generated by Thomas Kuhn's (1962) work on
the structure of scientific revolutions. Applying Popper's favorite
keyword of fallibilism, the statements of earlier scientists that
historiography had declared to be false were particularly attractive
objects for such an epistemologically guided revision.

The philosopher Imre Lakatos (1922-1972) was responsible for
introducing these new approaches into the history of mathematics.  One
of the examples he analyzed and published in 1966 received a great
deal of attention: Cauchy's theorem and the problem of uniform
convergence. Lakatos refines Robinson's approach by claiming that
Cauchy's theorem had also been correct at the time, because he had
been working with infinitesimals \cite[p.~431--432]{Sch}.
\end{quote}
One might have expected that, having devoted so much space to the
philosophical underpinnings of Lakatos' interpretation of Cauchy's sum
theorem, Schubring would actually devote a thought or two to that
interpretation itself.  Instead, Schubring presents a misguided claim
to the effect that Robinson acknowledged the incorrectness of the sum
theorem.%
\footnote{Schubring's quote of Robinson's reference to Cauchy's
``famous error'' is taken out of context.  Note that the full quote is
``a famous error of Cauchy's, which has been discussed repeatedly in
the literature'' (Robinson \cite[p.~271]{Ro66}).  Robinson devotes a
lengthy paragraph on pages 260-261 to a statement of a received view
of the history of the calculus.  Robinson then proceeds to refute the
received view.  Similarly, the received, and famous, view (and one
that ``has been discussed repeatedly in the literature'') is that
Cauchy erred in his treatment of the sum theorem.  Robinson proceeds
to challenge such a view, by offering an interpretation that
vindicates Cauchy's sum theorem, along the lines of Cauchy's own
modification/clarification of 1853 (Cauchy \cite{Ca53}).  Robinson's
interpretation is consistent with our position that Cauchy's 1821
result is irreducibly ambiguous, see also footnote~\ref{GG}.  For a
summary of the controversy over the sum theorem, see
Section~\ref{sum}.}
Schubring appears to feel that calling Lakatos a Popperian and a
fallibilist is sufficient refutation in its own right.  Similarly,
Schubring dismisses Laugwitz's reading of Cauchy as ``solipsistic"
\cite[p.~434]{Sch}; accuses them of interpreting Cauchy's conceptions
as
\begin{quote}
some hermetic closure of a {\em private\/} mathematics
\cite[p.~435]{Sch} [emphasis in the original--the authors];
\end{quote}
as well as being ``highly anomalous or isolated'' \cite[p.~441]{Sch}.
Now common sense would suggest that Laugwitz is interpreting Cauchy's
words according to their plain meaning, and takes his infinitesimals
at face value.  Isn't the burden of proof on Schubring to explain why
the triumvirate interpretation of Cauchy is not ``solipsistic'',
``hermetic'', or ``anomalous''?  Schubring does nothing of the sort.

Why are Lakatos and Laugwitz demonized rather than analyzed by
Schubring?  The issue of whether or not Schubring commands a minimum
background necessary to understand either Robinson's, Lakatos', or
Laugwitz's interpretation was discussed above.  More fundamentally,
the act of contemplating for a moment the idea that Cauchy's
infinitesimals can be taken at face value is unthinkable to a
triumvirate historian, as it would undermine the nominalistic
Cauchy-Weierstrass tale that the received historiography is erected
upon.  The failure to appreciate the potency of the
Robinson-Lakatos-Laugwitz interpretation is symptomatic of an ostrich
effect conditioned by a narrow A-continuum vision.%
\footnote{Such an effect is comparable to those occurring in the
constructivist context, see footnote~\ref{ostrich1}, and in the
educational context, see footnote~\ref{ostrich2}.}
The Robinson-Lakatos-Laugwitz interpretation of Cauchy's sum theorem
is considered in more detail in Section~\ref{sum}.

\section{A case study in nominalistic hermeneutics}

Chapter 6 in Schubring \cite{Sch} is entitled ``Cauchy's compromise
concept''.  Which compromise is the author referring to?  The answer
becomes apparent on page~439, where the author quotes Cauchy's
``reconciliation'' sentence:
\begin{quote}
My main aim has been to {\em reconcile\/} rigor, which I have made a
law in my Cours d'Analyse, with the simplicity that comes from the
direct consideration of infinitely small quantities (Cauchy 1823, see
\cite[p.~10]{Ca23}) [emphasis added--authors].
\end{quote}
Cauchy's choice of the word ``reconcile'' does suggest a resolution of
a certain tension.  What is the nature of such a tension?  The
sentence mentions ``rigor'' and ``infinitely small quantities'' in the
same breath.  This led Schubring to a conclusion of Cauchy's alleged
perception of a ``disagreement'' between them:
\begin{quote}
In his next textbook on differential calculus in 1823, Cauchy points
out expressly that he has adopted a compromise concept and that the
``simplicity of the infinitely small quantities'' [...]  disagrees
with the ``rigor" that he wished to achieve in his 1821 textbook
\cite[p.~439]{Sch}.
\end{quote}
Schubring's conclusion concerning such an alleged ``disagreement'', as
well as the ``compromise'' of his title, both hinge essentially on a
single word {\em concilier\/} (reconcile) in Cauchy.  Let us analyze
its meaning.  If it refers to a disagreement between rigor and
infinitesimals, how do we account for Cauchy's attribution, in 1821,
of a fundamental foundational role of infinitesimals in establishing a
rigorous basis for analysis?  Had Cauchy changed his mind sometime
between 1821 and 1823?

To solve the riddle we must place Cauchy's ``reconciliation'' sentence
in the context where it occurs.  In the sentence immediately preceding
it, Cauchy speaks of his break with the earlier texts in analysis:
\begin{quote}
Les m\'ethodes que j'ai suivies diff\`erent \`a plusieurs \'egards de
celles qui se trouvent expos\'ees dans les ouvrages du m\^eme genre
\cite[p.~10]{Ca23}.
\end{quote}
Could he be referring to his own earlier text?  To answer the
question, we must read on what Cauchy has to say.  Immediately
following the ``reconciliation'' sentence, Cauchy unleashes a
sustained attack against the flawed method of divergent power series.
Cauchy does not name the culprit, but clearly identifies the offending
treatise.  It is the {\em M\'ecanique Analytique\/}, see (Cauchy 1823,
p.~11).  The second edition of Lagrange's treatise came out in 1811,
when Cauchy was barely out of his teens.  Here Lagrange writes:
\begin{quote}
Lorsqu'on a bien con\c cu l'esprit de ce syst\`eme, et qu'on s'est
convaincu de l'exactitude de ses r\'esultats par la m\'ethode
g\'eom\'etrique des premi\`eres et derni\`eres raisons, ou par la
m\'ethode analytique des fonctions d\'eriv\'ees, on peut employer les
infiniment petits comme un instrument s\^ur et commode pour abr\'eger
et simplifier les d\'emonstrations \cite[p.~iv]{Lag}.
\end{quote}
Lagrange's ringing endorsement of infinitesimals in 1811 is as
unambivalent as that of Johann Bernoulli, l'Hopital, or Varignon.  In
rejecting Lagrange's flawed method of power series, as well as his
principle of the ``generality of algebra'', Cauchy was surely faced
with a dilemma with regard to Lagrange's infinitesimals, which had
stirred controversy for over a century.  We argue that it is the
context of a critical re-evaluation of Lagrange's mathematics that
created a tension for Cauchy vis-a-vis Lagrange's work of 1811: can he
sift the chaff from the grain?  

Cauchy's great accomplishment was his recognition that, while
Lagrange's flawed power series method and his principle of the
generality of algebra do not measure up to the standard of rigor
Cauchy sought to uphold in his own work, the infinitesimals can indeed
be reconciled with such a standard of rigor. The resolution of the
tension between the rejection of Lagrange's conceptual framework, on
the one hand, and the acceptance of his infinitesimals, on the other,
is what Cauchy is referring to in his ``reconciliation'' sentence.
Cauchy's blending of rigor and infinitesimals in 1823 is consistent
with his approach in 1821.  Cauchy's sentence compromises Schubring's
concept of a Cauchyan ambivalence with regard to infinitesimals, and
pulls the rug from under Schubring's nominalistic and solipsistic
reading of Cauchy.

\section{Cauchy's sum theorem}
\label{sum}

In this section, we summarize the controversy over the sum theorem,
recently analyzed by Br\aa ting \cite{Br}.  The issue hinges on two
types of convergence.  To clarify the mathematical issues involved, we
will first consider the simpler distinction between continuity and
uniform continuity.  Let~$x$ be in the domain of a function~$f$, and
consider the following condition, which we will call {\em
microcontinuity\/} at $x$:
\begin{quote}
``if~$x'$ is in the domain of~$f$ and~$x'$ is infinitely close to~$x$,
then~$f(x')$ is infinitely close to~$f(x)$".
\end{quote}
Then ordinary continuity of~$f$ is equivalent to~$f$ being
microcontinuous on the Archimedean continuum (A-continuum for short),
i.e., at every point $x$ of its domain in the A-continuum.  Meanwhile,
uniform continuity of~$f$ is equivalent to~$f$ being microcontinuous
on the Bernoulian continuum (B-continuum for short), i.e., at every
point $x$ of its domain in the B-continuum.%
\footnote{The relation of the two continua is discussed in more detail
in Section~\ref{rival1}.}

Thus, the function
\[
\sin \frac{1}{x}
\]
for positive~$x$ fails to be uniformly continuous because
microcontinuity fails at a positive infinitesimal~$x$.  The
function~$x^2$ fails to be uniformly continuous because of the failure
of microcontinuity at a single infinite member of the B-continuum.

A similar distinction exists between pointwise convergence and uniform
convergence.  The latter condition requires convergence at the points
of the B-continuum in addition to the points of the A-continuum, see
e.g.~Goldblatt \cite[Theorem 7.12.2, p.~87]{Go}.  

Which condition did Cauchy have in mind in 1821?  This is essentially
the subject of the controversy over the sum theorem.

Abel interpreted it as convergence on the A-continuum, and presented
``exceptions'' (what we would call today counterexamples) in 1826.
After the publication of additional such exceptions by Seidel and
Stokes in the 1840s, Cauchy clarified/modified his position in 1853.
In his text \cite{Ca53}, he specified a stronger condition of
convergence on the B-continuum, including at~$x=1/n$.  The latter
entity is explicitly mentioned by Cauchy as illustrating the failure
of the error term to tend to zero.  The stronger condition bars Abel's
counterexample.  See our text \cite{KK11b} for more details.
%
%

\section{Fermat, Wallis, and an ``amazingly reckless'' use of infinity}
\label{rival1}

\begin{figure}
\begin{equation*}
\xymatrix@C=95pt{{} \ar@{-}[rr] \ar@{-}@<-0.5pt>[rr]
\ar@{-}@<0.5pt>[rr] & {} \ar@{->}[d]^{\hbox{st}} & \hbox{\quad
B-continuum} \\ {} \ar@{-}[rr] & {} & \hbox{\quad A-continuum} }
\end{equation*}
\caption{\textsf{Thick-to-thin: taking standard part (the thickness of
the top line is merely conventional, and meant to suggest the presence
of additional numbers, such as infinitesimals)}}
\label{31}
\end{figure}

A Leibnizian definition of the derivative as the infinitesimal
quotient
\begin{equation*}
\frac{\Delta y}{\Delta x},
\end{equation*} 
whose logical weakness was criticized by Berkeley, was modified by
A.~Robinson by exploiting a map called {\em the standard part\/},
denoted~``st'', from the finite part of a B-continuum (for
``Bernoullian''), to the A-continuum (for ``Archimedean''), as
illustrated in Figure~\ref{31}.%
\footnote{In the context of the hyperreal extension of the real
numbers, the map ``st'' sends each finite point~$x$ to the real point
st$(x)\in \R$ infinitely close to~$x$.  In other words, the map ``st''
collapses the cluster (halo) of points infinitely close to a real
number~$x$, back to~$x$.}
Here two points of a B-continuum have the same image under ``st'' if
and only if they are equal up to an infinitesimal.  

This section analyzes the historical seeds of Robinson's theory, in
the work of Fermat, Wallis, as well as Barrow.%
\footnote{\label{barrow}While Barrow's role is also critical, we will
mostly concentrate on Fermat and Wallis.}
The key concept here is that of {\em adequality\/} (see below).  It
should be kept in mind that Fermat never considered the local slope of
a curve.  Therefore one has to be careful not to attribute to Fermat
mathematical content that could not be there.  On the other hand,
Barrow did study curves and their slope.  Furthermore, Barrow
exploited Fermat's adequality in his work~\cite[p.~252]{Barr}, as
documented by H.~Breger \cite[p.~198]{Bre94}.

The binary relation of ``equality up to an infinitesimal'' was
anticipated in the work of Pierre de Fermat.  Fermat used a term
usually translated into English as ``adequality''.%
\footnote{In French one uses {\em ad\'egalit\'e, ad\'egal\/}, see
\cite[p.~73]{Je}.}
Andr\'e Weil writes as follows:
\begin{quote} 
Fermat [...] developed a method which slowly but surely brought him
very close to modern infinitesimal concepts.  What he did was to write
congruences between functions of~$x$ modulo suitable powers
of~$x-x_0$; for such congruences, he introduces the technical term
{\em adaequalitas, adaequare\/}, etc., which he says he has borrowed
from Diophantus.  As Diophantus V.11 shows, it means an approximate
equality, and this is indeed how Fermat explains the word in one of
his later writings \cite[p.~1146]{We}.
\end{quote} 
Weil \cite[p.~1146, footnote~5]{We} then supplies the following quote
from Fermat:
\begin{quote}
{\em Adaequetur, ut ait Diophantus,%
\footnote{The original term in Diophantus is~$\pi \alpha \rho \iota
\sigma \acute{o} \tau \eta \varsigma$, see Weil \cite[p.~28]{We84}.}
aut fere aequetur\/}; in Mr.~Mahoney's translation: ``adequal, or
almost equal" (p. 246).
\end{quote}
Here Weil is citing Mahoney \cite[p.~246]{Mah73}
(cf. \cite[p.~247]{Mah94}).  Mahoney similarly mentions the meaning of
``approximate equality'' or ``equality in the limiting case'' in
\cite[p.~164, end of footnote~46]{Mah73}.  Mahoney also points out
that the term ``adequality'' in Fermat has additional meanings.  The
latter are emphasized in a recent text by E.~Giusti~\cite{Giu}, who is
sharply critical of Breger~\cite{Bre94}.  While the review~\cite{We}
by Weil is similarly sharply critical of Mahoney, both agree that the
meaning of ``approximate equality'', leading into infinitesimal
calculus, is at least {\em one of the meanings\/} of the term {\em
adequality\/} for Fermat.%
\footnote{Jensen similarly describes adequality as approximate
equality, and describes neglected terms as {\em infinitesimals\/} in
\cite[p.~82]{Je}.  Struik notes that ``Fermat uses the term to denote
what we call a limiting process'' \cite[p.~220, footnote~5]{Stru}.
K.~Barner~\cite{Ba} compiled a useful bibliography on Fermat's
adequality, including many authors we have not mentioned here.}

\begin{figure}
\includegraphics[height=2.3in]{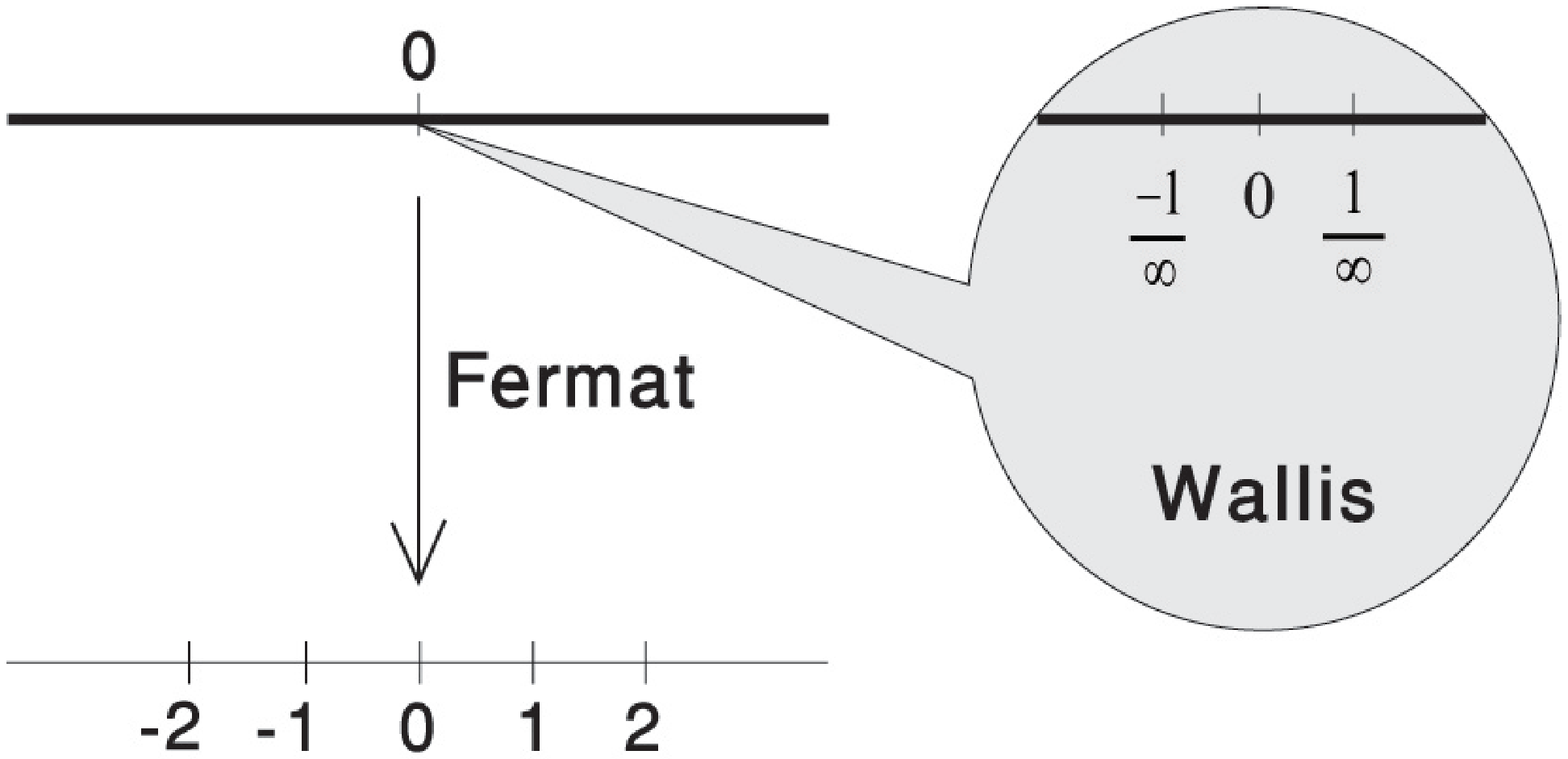}
\caption{\textsf{Zooming in on Wallis's
infinitesimal~$\frac{1}{\infty}$, which is adequal to~$0$ in Fermat's
terminology}}
\label{FermatWallis}
\end{figure}

This meaning was aptly summarized by J.~Stillwell.  Stillwell's
historical presentation is somewhat simplified, and does not
sufficiently distinguish between the seeds actually present in Fermat,
on the one hand, and a modern interpretation thereof, on the other,%
\footnote{See main text around footnote~\ref{barrow} above for a
discussion of Barrow's role, documented by H~Breger.}
%
%
but he does a splendid job of explaining the mathematical background
for the uninitiated.  Thus, he notes that~$2x+dx$ is not equal to~$2x$
(see Figure~\ref{jul10}), and writes:
\begin{quote}
Instead, the two are connected by a looser notion than equality that
Fermat called adequality.  If we denote adequality by~$=_{ad}$, then
it is accurate to say that
\[
2x+dx=_{ad}2x,
\]
and hence that~$dy/dx$ for the parabola is adequal to~$2x$.
Meanwhile,~$2x+dx$ is not a number, so~$2x$ is the only number to
which~$dy/dx$ is adequal.  This is the true sense in which~$dy/dx$
represents the slope of the curve~\cite[p.~91]{Sti}.
\end{quote}
Stillwell points out that
\begin{quote}
Fermat introduced the idea of adequality in 1630s but he was ahead of
his time.  His successors were unwilling to give up the convenience of
ordinary equations, preferring to use equality loosely rather than to
use adequality accurately.  The idea of adequality was revived only in
the twentieth century, in the so-called non-standard analysis
\cite[p.~91]{Sti}.
\end{quote}
We will refer to the map from the (finite part of the) B-continuum to
the A-continuum as the Fermat-Robinson standard part, see
Figure~\ref{FermatWallis}.

As far as the logical criticism formulated by Rev.~George is
concerned, Fermat's adequality had pre-emptively provided the seeds of
an answer, a century before the bishop ever lifted up his pen to write
{\em The Analyst\/}~\cite{Be}.

\begin{figure}
\includegraphics[height=2.3in]{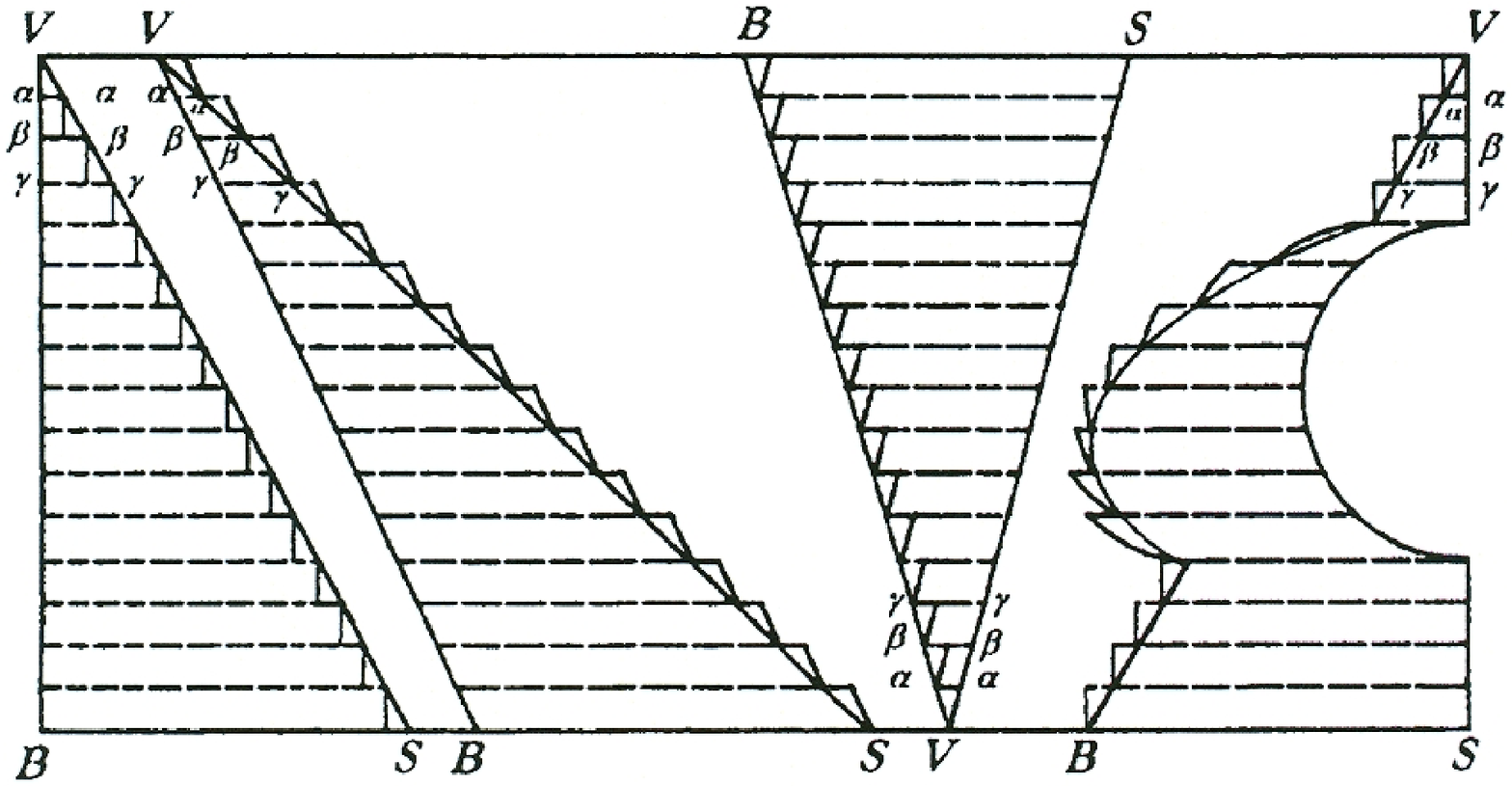}
\caption{\textsf{Area calculations in Wallis: slicing it up into
dilatable parallelograms of infinitesimal altitude}}
\label{Wallis}
\end{figure}

Fermat's contemporary John Wallis, in a departure from Cavalieri's
focus on the geometry of indivisibles, emphasized the arithmetic of
infinitesimals, see J.~Stedall's introduction in \cite{Wa}.  To
Cavalieri, a plane figure is made of lines; to Wallis, it is made of
parallelograms of infinitesimal altitude.  Wallis transforms this
insight into symbolic algebra over the~$\infty$ symbol which he
introduced.  He exploits formulas
like~$\infty\times\frac{1}{\infty}=1$ in his calculations of areas.
Thus, in proposition 182 of {\em Arithmetica Infinitorum\/}, Wallis
partitions a triangle of altitude~$A$ and base~$B$ into a precise
number~$\infty$ of ``parallelograms" of infinitesimal
width~$\frac{A}{\infty}$, see Figure~\ref{Wallis} (copied from
\cite[p.~170]{Mun}).

He then computes the combined length of the bases of the
parallelograms to be~$\frac{B}{2}\infty$, and finds the area to be
\begin{equation}
\label{wallis}
\frac{A}{\infty}\times \frac{B}{2}\infty = \frac{AB}{2}.
\end{equation}
Wallis used an actual infinitesimal~$\frac{1}{\infty}$ in
calculations as if it were an ordinary number, anticipating Leibniz's
law of continuity.

Wallis's area calculation~\eqref{wallis} is reproduced by J.~Scott,
who notes that Wallis
\begin{quote}
treats infinity as though the ordinary rules of arithmetic could be
applied to it \cite[p.~20]{Sco}.
\end{quote}
Such a treatment of infinity strikes Scott as something of a blemish,
as he writes:
\begin{quote}
But this is perhaps understandable.  For many years to come the
greatest confusion regarding these terms persisted, and even in the
next century they continued to be used in what appears to us an
amazingly reckless fashion \cite[p.~21]{Sco}.
\end{quote}
What is the source of Scott's confidence in dismissing Wallis's use of
infinity as ``reckless''?  Scott identifies it on the preceding page
of his book; it is, predictably, the triumvirate ``modern conception
of infinity'' \cite[p.~19]{Sco}.  Scott's tunnel A-continuum vision
blinds him to the potential of Wallis's vision of infinity.  But this
is perhaps understandable.  Many years separate Scott from Robinson's
theory which in particular empowers Wallis's calculation.  The lesson
of Scott's condescending steamrolling of Wallis's infinitesimal
calculation could be taken to heart by historians who until this day
cling to a nominalistic belief that Robinson's theory has little
relevance to the history of mathematics in the 17th century.

\section{Conclusion}

Nominalism in the narrow sense defines its ontological target as the
ordinary numbers.  Burgess in his essay \cite{Bu83} suggests that
there is also a nominalism in a broader sense.  Thus, he quotes at
length Manin's criticism of constructivism, suggesting that
LEM-elimination can also fall under the category of a nominalism
understood in a broader sense.  In a later text \cite[p.~30]{Bu04},
Burgess discusses Brouwer under a similar angle.

Infinitesimals were largely eliminated from mathematical discourse
starting in the 1870s through the efforts of the great triumvirate.%
\footnote{See footnote~\ref{triumvirate}.}
The elimination took place under the banner of striving for greater
rigor,%
\footnote{See Vicenti and Bloor \cite{VB} for an analysis of rigor in
19th century mathematics.  See also our Section~\ref{six}.}
but the roots of the triumvirate reconstruction lay in a failure to
provide a solid foundation for a B-continuum (see
Section~\ref{rival1}).  Had actually useful mathematics been
sacrificed on the altar of ``mathematical rigor" during the second
half of the 19th century?

Today we can give a precise sense to C.S. Peirce's description of the
real line as a {\em pseudo-continuum\/}.%
\footnote{American philosopher Charles Sanders Peirce felt that a
construction of a true continuum necessarily involves infinitesimals.
He wrote as follows: ``But I now define a pseudo-continuum as that
which modern writers on the theory of functions call a continuum.  But
this is fully represented by [...]  the totality of real values,
rational and irrational'' \cite{Pe} (see CP 6.176, 1903 marginal note.
Here, and below, CP x.y stands for Collected Papers of Charles Sanders
Peirce, volume x, paragraph y).  Peirce used the word
``pseudo-continua" to describe real numbers in the syllabus (CP 1.185)
of his lectures on Topics of Logic.  Thus, Peirce's intuition of the
continuum corresponded to a type of a B-continuum (see
Section~\ref{rival1}), whereas an A-continuum to him was a
pseudo-continuum.}

Cantor's revolutionary advances in set theory went hand-in-hand with
his emotional opposition to infinitesimals as an ``abomination'' and
the ``cholera bacillus'' of mathematics.  Cantor's interest in
eliminating infinitesimals is paralleled nearly a century later by
Bishop's interest in eliminating LEM, and by a traditional
nominalist's interest in eliminating Platonic counting numbers.  The
automatic infinitesimal-to-limit translation as applied to Cauchy by
Boyer and others is not only reductionist, but also
self-contradictory, see \cite{KK11b}.


\section{Rival continua}
\label{rival2}

This section summarizes a 20th century implementation of the
B-continuum, not to be confused with incipient notions of such a
continuum found in earlier centuries.  An alternative implementation
has been pursued by Lawvere, John L. Bell \cite{Bel08, Bel}, and
others.

We illustrate the construction by means of an infinite-resolution
microscope in Figure~\ref{FermatWallis}.  We will denote such a
B-continuum by the new symbol \RRR{} (``thick-R'').  Such a continuum
is constructed in formula~\eqref{bee}.  We will also denote its finite
part, by
\begin{equation*}
\RRR_{<\infty} = \left\{ x\in \RRR : \; |x|<\infty \right\},
\end{equation*}
so that we have a disjoint union
\begin{equation}
\RRR= \RRR_{<\infty} \cup \RRR_{\infty},
\end{equation}
where~$\RRR_{\infty}$ consists of unlimited hyperreals (i.e., inverses
of nonzero infinitesimals).

The map ``st'' sends each finite point~$x\in \RRR$, to the real point
st$(x)\in \R$ infinitely close to~$x$, as follows:%
\footnote{This is the Fermat-Robinson standard part whose seeds in
Fermat's adequality were discussed in Appendix~\ref{rival1}.}
%
%
\begin{equation*}
\xymatrix{\quad \RRR_{{<\infty}}^{~} \ar[d]^{{\rm st}} \\ \R}
\end{equation*}
Robinson's answer to Berkeley's {\em logical criticism\/} (see
D.~Sherry \cite{She87}) is to define the derivative as
\begin{equation*}
\hbox{st} \left( \frac{\Delta y}{\Delta x} \right),
\end{equation*}
instead of~$\Delta y/\Delta x$.

Note that both the term ``hyper-real field'', and an ultrapower
construction thereof, are due to E.~Hewitt in 1948, see
\cite[p.~74]{Hew}.  In 1966, Robinson referred to the 
\begin{quote}
theory of hyperreal fields (Hewitt [1948]) which ... can serve as
non-standard models of analysis \cite[p.~278]{Ro66}.
\end{quote}
The {\em transfer principle\/} is a precise implementation of
Leibniz's heuristic {\em law of continuity\/}: ``what succeeds for the
finite numbers succeeds also for the infinite numbers and vice
versa'', see~\cite[p.~266]{Ro66}.  The transfer principle, allowing an
extention of every first-order real statement to the hyperreals, is a
consequence of the theorem of J.~{\L}o{\'s} in 1955, see~\cite{Lo},
and can therefore be referred to as a Leibniz-{\L}o{\'s} transfer
principle.  A Hewitt-{\L}o{\'s} framework allows one to work in a
B-continuum satisfying the transfer principle.  To elaborate on the
ultrapower construction of the hyperreals, let~$\Q^\N$ denote the ring
of sequences of rational numbers.  Let
\begin{equation*}
\left( \Q^\N \right)_C
\end{equation*}
denote the subspace consisting of Cauchy sequences.  The reals are by
definition the quotient field
\begin{equation}
\label{real}
\R:= \left. \left( \Q^\N \right)_C \right/ \mathcal{F}_{\!n\!u\!l\!l},
\end{equation}
where~$\mathcal{F}_{\!n\!u\!l\!l}$ contains all null sequences.
Meanwhile, an infinitesimal-enriched field extension of~$\Q$ may be
obtained by forming the quotient
\begin{equation*}
\left.  \Q^\N \right/ \mathcal{F}_{u}.
\end{equation*}
Here a sequence~$\langle u_n : n\in \N \rangle$ is
in~$\mathcal{F}_{u}$ if and only if the set of indices
\[
\{ n \in \N : u_n = 0 \}
\]
is a member of a fixed ultrafilter.%
\footnote{In this construction, every null sequence defines an
infinitesimal, but the converse is not necessarily true.  Modulo
suitable foundational material, one can ensure that every
infinitesimal is represented by a null sequence; an appropriate
ultrafilter (called a {\em P-point\/}) will exist if one assumes the
continuum hypothesis, or even the weaker Martin's axiom.  See Cutland
{\em et al\/} \cite{CKKR} for details.}
See Figure~\ref{helpful}.

\begin{figure}
\begin{equation*}
\xymatrix{ && \left( \left. \Q^\N \right/ \mathcal{F}_{\!u}
\right)_{<\infty} \ar@{^{(}->} [rr]^{} \ar@{->>}[d]^{\rm st} &&
\RRR_{<\infty} \ar@{->>}[d]^{\rm st} \\ \Q \ar[rr] \ar@{^{(}->} [urr]
&& \R \ar[rr]^{\simeq} && \R }
\end{equation*}
\caption{\textsf{An intermediate field~$\left. \Q^\N \right/
\mathcal{F}_{\!u}$ is built directly out of~$\Q$}}
\label{helpful}
\end{figure}

To give an example, the sequence
\begin{equation}
\label{infinitesimal}
\left\langle \tfrac{(-1)^n}{n} \right\rangle
\end{equation}
represents a nonzero infinitesimal, whose sign depends on whether or
not the set~$2\N$ is a member of the ultrafilter.  To obtain a full
hyperreal field, we replace~$\Q$ by~$\R$ in the construction, and form
a similar quotient
\begin{equation}
\label{bee}
\RRR:= \left.  \R^\N \right/ \mathcal{F}_{u}.
\end{equation}
We wish to emphasize the analogy with formula~\eqref{real} defining
the A-continuum.  Note that, while the leftmost vertical arrow in
Figure~\ref{helpful} is surjective, we have
\begin{equation*}
\left( \Q^\N / \mathcal{F}_{u} \right) \cap \R = \Q.
\end{equation*}
A more detailed discussion of this construction can be found in the
book by M.~Davis~\cite{Dav77}.  
%
%
See also B\l aszczyk \cite{Bl} for some philosophical implications.
More advanced properties of the hyperreals such as saturation were
proved later, see Keisler \cite{Ke94} for a historical outline.  A
helpful ``semicolon'' notation for presenting an extended decimal
expansion of a hyperreal was described by A.~H.~Lightstone~\cite{Li}.
See also P.~Roquette \cite{Roq} for infinitesimal reminiscences.  A
discussion of infinitesimal optics is in K.~Stroyan \cite{Str},
J.~Keisler~\cite{Ke}, D.~Tall~\cite{Ta80}, and L.~Magnani and
R.~Dossena~\cite{MD, DM}, and Bair \& Henry \cite{BH}.  Applications
of the B-continuum range from aid in teaching calculus \cite{El, KK1,
KK2, Ta91, Ta09a} (see illustration in Figure~\ref{jul10}) to the
Bolzmann equation (see L.~Arkeryd~\cite{Ar81, Ar05}); modeling of
timed systems in computer science (see H.~Rust \cite{Rust});
mathematical economics (see Anderson \cite{An00}); mathematical
physics (see Albeverio {\em et al.\/} \cite{Alb}); etc.

\section*{Acknowledgments}

We are grateful to Martin Davis, Solomon Feferman, Reuben Hersh, David
Sherry, and Steve Shnider for invaluable comments that helped improve
the manuscript.  Hilton Kramer's influence is obvious throughout.

\end{document}